# Computation of a universal deformation ring in the Borel case

By ARIANE MÉZARD

*Institut Fourier, Université de Grenoble 1*

*38402 Saint-Martin d'Hères, France*

*e-mail:* `mezard@ujf-grenoble.fr`

*Abstract*

We compute the universal deformation ring of an odd Galois representation $\bar{\rho} : \mathrm{Gal}(M/\mathbf{Q}) \to \mathrm{Gl}_2(\mathbf{F}_p)$ with an upper triangular image, where $M$ is the maximal abelian pro-$p$-extension of $F_\infty$ unramified outside a finite set of places $S$, $F_\infty$ being a free pro-$p$-extension of a subextension $F$ of the field $K$ fixed by $\mathrm{Ker}\bar{\rho}$. We establish a link between the latter universal deformation ring and the universal deformation ring of $\bar{\rho} : \mathrm{Gal}(K_S/\mathbf{Q}) \to \mathrm{Gl}_2(\mathbf{F}_p)$, where $K_S$ is the maximal pro-$p$-extension of $K$ unramified outside $S$. We then give some examples.

*Résumé*

On calcule l'anneau de déformation universel d'une représentation galoisienne impaire $\bar{\rho} : \mathrm{Gal}(M/\mathbf{Q}) \to \mathrm{Gl}_2(\mathbf{F}_p)$ d'image triangulaire supérieure où $M$ est la plus grande pro-$p$-extension abélienne non ramifiée en dehors d'un ensemble fini $S$ de places de $F_\infty$, pour $F_\infty$ une pro-$p$-extension d'un sous corps $F$ du corps $K$ fixe par le noyau de $\bar{\rho}$. On établit un lien entre cet anneau de déformation universel et l'anneau de déformation universel de $\bar{\rho} : \mathrm{Gal}(K_S/\mathbf{Q}) \to \mathrm{Gl}_2(\mathbf{F}_p)$, pour $K_S$ la plus grande pro-$p$-extension de $K$ non ramifiée en dehors de $S$ ; et on donne quelques exemples.

## 1. *Introduction*

Let $p$ be an odd prime number. Let $\bar{\rho} : \mathrm{Gal}(\bar{\mathbf{Q}}/\mathbf{Q}) \to \mathrm{Gl}_2(\mathbf{F}_p)$ be an odd continuous representation unramified outside a finite set $S_{\mathbf{Q}}$ of rational primes. In that situation $\bar{\rho}$ factors through $G_S = \mathrm{Gal}(K_{S,p}/\mathbf{Q})$ where $K$ denotes the field fixed by $\mathrm{Ker}\bar{\rho}$ in $\bar{\mathbf{Q}}$, $S$ a finite set of places of $K$ containing the places over $S_{\mathbf{Q}}$ in $K$ and $K_{S,p}$ the maximal pro-$p$-extension of $K$, unramified outside $S$.

Mazur's deformation theory of Galois representations shows the existence of a (uni)versal deformation ring $R_{G_S}(\bar{\rho})$ and the associated (uni)versal representation $\rho_{G_S}$ which allow us to parametrize all deformations $\rho : G_S \to \mathrm{Gl}_2(R)$ of $\bar{\rho} : G_S \to \mathrm{Gl}_2(\mathbf{F}_p)$ where $R$ stands for any complete noetherian local ring with residual field $\mathbf{F}_p$.

To our knowledge, there is not yet any general method for computing the (uni)versal deformation ring $R_{G_S}(\bar{\rho})$ although partial results suggest that the structure of this ring is closely related to the adjoint representation $Ad(\bar{\rho})$. By application of Schlessinger's criterion [Sc] Mazur ([Ma] subsection 1.2) shows that $R_{G_S}(\bar{\rho})$ is a quotient ring of a formal series ring whose minimal number of variables is $d' = \dim_{\mathbf{F}_p} H^1(G_S, Ad(\bar{\rho}))$. To be more specific, $R_{G_S}(\bar{\rho}) = \mathbf{Z}_p[[Y_1, \ldots, Y_{d'}]]/I$ where $I$ will be called the *ideal of relations* of $R_{G_S}(\bar{\rho})$, so that the determination of $R_{G_S}(\bar{\rho})$ amounts to that of $I$.

The knowledge of $I$ would allow us to discuss Mazur's question on the Krull dimension of $R_{G_S}(\bar{\rho})/pR_{G_S}(\bar{\rho})$ (in the case dim$\bar{\rho} = 1$, this is the celebrated Leopoldt conjecture, [Ma] subsections 1.6 and 1.10). The general theory of obstructions shows that $I = 0$ (i.e the (uni)versal deformation ring is free) if $H^2(G_S, Ad(\bar{\rho})) = 0$. Almost all explicit examples that we know of rely on the latter assumption ([Bo1],[Ra]), which occurs for instance if the maximal pro-$p$-quotient $P_{G_S}$ of $G_S$ is a free pro-$p$-group or if some irreducible components of the representation $P_{G_S}^{ab} \otimes \mathbf{F}_p$ are prime to the irreducible components of the representation $Ad(\bar{\rho})$. This suggests a precise connection (via $\bar{\rho}$) between the relations of the pro-$p$-group $P_{G_S}$ and the relations of the (uni)versal deformation ring $R_{G_S}$. It should be stressed though, that even the precise knowledge of the relations of $P_{G_S}$ does not imply straightforwardly the knowledge of the ideal $I$.

In this paper, we present an approach for determining $I$, via Iwasawa theory, built on an example studied by Boston. In [Bo1] subsection 9.3 Boston considers representations $\bar{\rho}$ which factor through $\bar{\rho}_G : G \to \mathrm{Gl}_2(\mathbf{F}_p)$, where $G$ is a natural quotient of $G_S$ which occurs in Iwasawa theory. By construction there exists a surjective morphism of local rings $R_{G_S}(\bar{\rho}) \to R_G(\bar{\rho}_G)$ making $R_G(\bar{\rho}_G)$ an approximation of $R_{G_S}(\bar{\rho})$. We make more precise this morphism in the sequel. Moreover if the order of Im$\bar{\rho}$ is prime to $p$ (this is the *tame case*), then $R_{G_S}(\bar{\rho})$ and $R_G(\bar{\rho}_G)$ have the same minimal number of variables.

In the sequel we shall develop Boston's example along two directions: First we present a systematic way of deriving, from a minimal system of relations of the group $P_{G_S}$, a minimal set of relations for the ring $R_G(\bar{\rho}_G)$, using the method of Fox derivatives described in [Ng1]. Second we extend the framework of Iwasawa theory to study not only the $\mathbf{Z}_p$-cyclotomic extension, but also pro-$p$-free extensions, which allows us to enlarge the group $G$ which approximates $G_S$ and to replace the classical Iwasawa algebra by the (non commutative) Magnus algebra as in [Ng2].

More precisely it will be shown that the computation of $R_G(\bar\rho_G)$ is possible provided that $\mathrm{Im}\bar\rho \subset \begin{pmatrix} * & * \\ 0 & * \end{pmatrix}$ with an additional hypothesis on the diagonal characters (this is the *Borel case*).

If the relations of $P_{G_S}$ are known, [Bo1] proposition 6.1 allows us in principle to determine $R_{G_S}(\bar\rho)$ in the tame case, though in most cases the computations involved cannot be carried out. Our method here works fine even if the relations of $P_{G_S}$ are not explicitly known.

The method involves two steps: first (sections 4-5) we find the images of the generators of $G$ using the action of $G/P_G$ (where $P_G$ is the $p$-Sylow subgroup of $G$). We proceed exactly as if $P_G$ were free. Second (sections 6-7) we express the relations between selected generators (the images of which are simple).

In the special case considered by Boston, we obtain a presentation of the universal ring $R_G(\bar\rho_G) = \mathbf{Z}_p[[Y_1, \ldots, Y_{d'_G}]]/I$ where $d'_G = \dim_{\mathbf{F}_p} H^1(G, Ad\bar\rho_G)$ and $I$ is the ideal of relations generated by

$$\Big(\prod_{i=1}^{u_{X_\infty}} (1+Y_i)^{a_i^j} - 1\Big)_{1\leq j\leq \ell}$$

and

$$\Big(a_n^j + \sum_{k=1}^{\infty} b_{n,k}^j \Big(\frac{1+Y_{d'_G-1}}{1+Y_{d'_G}} - 1\Big)^k + \sum_{i=1+u_{X_\infty}}^{v_{X_\infty}} \Big(a_i^j + \sum_{k=1}^{\infty} b_{i,k}^j \Big(\frac{1+Y_{d'_G-1}}{1+Y_{d'_G}} - 1\Big)^k\Big) Y_i\Big)_{1\leq j\leq \ell}$$

where the coefficients $a_i^j, b_{i,k}^j$ are derived from the relations of $P_{G_S}$ (see theorem 7.2.1 and compare [Bo1] subsection 9.3).

The comparison between the universal deformation rings $R_G(\bar\rho_G)$ and $R_{G_S}(\bar\rho)$ is done in theorem 7.4.2, where we make precise the natural surjective morphism $R_{G_S}(\bar\rho) \to R_G(\bar\rho_G)$.

Lastly (sections 7-8) we give some applications of our results to representations associated to elliptic curves or to modular forms. For the latter representations, cyclotomic fields (we correct an imprecise result of [Bo1] proposition 9.2 and [Bo2] section 6) or fields having a so-called Wingberg presentation appear naturally.

## 2. *Notations*

Let $\bar{\rho} : \mathrm{Gal}(\bar{\mathbf{Q}}/\mathbf{Q}) \to \mathrm{Gl}_2(\mathbf{F}_p)$ be an odd continuous representation unramified outside a finite set of primes $S_{\mathbf{Q}}$ of $\mathbf{Q}$.
Assume that $\mathrm{Im}\bar{\rho}$ is upper triangular: $\mathrm{Im}\bar{\rho} \subset \begin{pmatrix} * & * \\ 0 & * \end{pmatrix}$.
Let $K$ denote the subfield of $\bar{\mathbf{Q}}$ fixed by $\mathrm{Ker}\bar{\rho}$. Hence, one has $\mathrm{Gal}(K/\mathbf{Q}) \cong \mathrm{Im}\bar{\rho}$ (finite). Let $F$ be a subextension of $K$ such that $\mathrm{Gal}(K/F)$ (possibly trivial) is the Sylow $p$-subgroup of $\mathrm{Gal}(K/\mathbf{Q})$.
We assume the extension $K/F$ is ramified at $p$.

Let $S_F$ be a finite set of places of $F$ containing the places over $S_{\mathbf{Q}}$ in $F$. Suppose $S_F$ contains the primes over $p$ in $F$ and the archimedean primes. Let $S_K$ be the places in $K$ over $S_F$. We write $S$ for $S_F$ or $S_K$ when no confusion can arise.

The maximal algebraic pro-$p$-extensions of $F$ or $K$ unramified outside $S_F$ or $S_K$ are the same: $K_{S,p} = F_{S,p}$. We work with the following setup

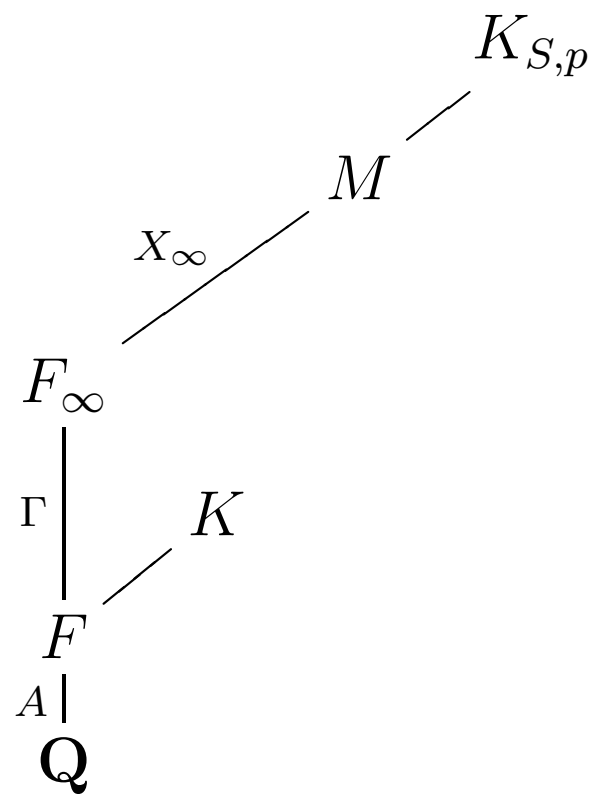

where $F_\infty$ is a free pro-$p$-extension (not necessary abelian) of rank $k \geq 1$ of $F$: this means that $\Gamma = \mathrm{Gal}(F_\infty/F)$ is pro-$p$-free on $k$-generators. We assume that $A$ acts on $\Gamma$ by conjugation.

The field $M$ is the maximal abelian pro-$p$-extension of $F_\infty$ unramified outside $S$,
and $H = \mathrm{Gal}(K_{S,p}/F_\infty)$,
and $X_\infty = \mathrm{Gal}(M/F_\infty) = H^{ab}$.

**Remark 2.1** We assume that $K$ and $F_\infty$ are linearly disjoint.

Define furthermore $G_S = \mathrm{Gal}(K_{S,p}/\mathbf{Q})$,
$P_{G_S} = \mathrm{Gal}(K_{S,p}/F)$,
$G = \mathrm{Gal}(M/\mathbf{Q})$,
$P_G = \mathrm{Gal}(M/F)$;
The group $\bar{P}_{G_S}$ (resp. $\bar{P}_G$, resp. $\bar{\Gamma}$) is the quotient of $P_{G_S}$ (resp. $P_G$, resp. $\Gamma$) by its Frattini subgroup.

**Remark 2.2** Even if $G$ and $P_G$ have no index $S$, they depend on $S_F$ and on the choice of $\Gamma$.

In the sequel we will work in particular with
-the $\mathbf{Z}_p$-cyclotomic extension of $F$, (then we denote $\Gamma = \Gamma_{cyc}$) in order to maximize information on the module $X_\infty$,
-a maximal free pro-$p$-extension of $F$, (then we denote $\Gamma = \Gamma_{max}$) in order to maximize the size of $G$. Note that such a maximal free pro-$p$-extension is not necessary unique (see [Ya2] remark).
In these cases it is clear that $A$ acts on $\Gamma$.

According to Schur-Zassenhaus theorem (see proposition 2.1 [Bo1]) the group $G$ contains a subgroup $A$ mapping isomorphically to $G/P_G$.
We assume

$$\bar{\rho}_{|A} = \begin{pmatrix} \chi_1 & 0 \\ 0 & \chi_2 \end{pmatrix} \text{ such that } \chi_1 \neq \pm\chi_2$$

where $\chi_1, \chi_2 : A \to \mathbf{F}_p^*$ are the diagonal characters defined by $\mathrm{Im}\bar{\rho}$. The case where $\mathrm{Im}\bar{\rho}$ is upper triangular and $\chi_1 \neq \chi_2$ is termed the *Borel case*.

We denote by $Ad(\bar{\rho})$ the vector space $\mathrm{M}_2(\mathbf{F}_p)$ equipped with the action of $G_S$ defined by conjugation by $\bar{\rho}$ (or the action of $G$). Depending upon the context we shall instead of $Ad(\bar{\rho})$ write $Ad(\bar{\rho}_G)$ when necessary, to avoid confusion.

At last we define $\Gamma_2(R) = \mathrm{Ker}(\mathrm{Gl}_2(R) \to \mathrm{Gl}_2(\mathbf{F}_p))$ for all $R \in \mathcal{C}$ where $\mathcal{C}$ is the category of the complete noetherian local rings of residual field $\mathbf{F}_p$. We denote by $\mathcal{M}_R$ the maximal ideal of $R$.

Two liftings $\rho$ and $\rho'$ to $R$ of $\bar{\rho}$ are *strictly equivalent* if there exists $M' \in \Gamma_2(R)$ such that $\rho' = M'\rho M'^{-1}$.

Mazur's functor of deformations of $\bar{\rho}$, $\mathcal{F} : \mathcal{C} \to \mathrm{Sets}$ defined by $\mathcal{F}(R)$ is the set of the strict equivalent classes of liftings of $\bar{\rho}$ onto $R$. One takes $\rho \in \mathcal{F}(R)$.

There exists a canonical surjective morphism $G_S \to G$ and, by construction, the representation $\bar{\rho}$ factors through $G$; after [Ma] subsection 1.3 there exists a surjective morphism of local rings $R_{G_S}(\bar{\rho}) \to R_G(\bar{\rho}_G)$.

**Remark 2.3** We have the exact sequence:

$$0 \to P_{G_S} \to G_S \to \mathrm{Gal}(F/\mathbf{Q}) \to 0$$

In the tame case, by definition, the order of $\mathrm{Gal}(F/\mathbf{Q})$ is prime to $p$ then $H^1(G_S, Ad(\bar{\rho})) \cong H^1(P_{G_S}, Ad(\bar{\rho}))^{G_S/P_{G_S}}$. Moreover in the tame case, $P_{G_S} \subset \mathrm{Ker}\bar{\rho}$, i.e the action of $P_{G_S}$ on $Ad(\bar{\rho})$ is trivial; hence

$$H^1(P_{G_S}, Ad(\bar{\rho}))^{G_S/P_{G_S}} \cong \mathrm{Hom}(\bar{P}_{G_S}, Ad(\bar{\rho}))^{G_S/P_{G_S}} \cong H^1(G_S, Ad(\bar{\rho}))$$

This reasoning is also valuable if we replace $G_S$ and $Ad(\bar{\rho})$ by $G$ resp. $Ad(\bar{\rho}_G)$, and since $\bar{P}_{G_S} \cong \bar{P}_G$ we obtain

$$H^1(G_S, Ad(\bar{\rho})) \cong H^1(G, Ad(\bar{\rho}_G))$$

Hence $R_{G_S}(\bar{\rho})$ and $R_G(\bar{\rho}_G)$ have the same minimal number $d'$ of variables in the tame case (for example if $\mathrm{Im}\bar{\rho}$ is diagonal).

Recall finally that, following [Ra] theorem 1.1, if $Ad(\bar{\rho})^{G_S} = \mathbf{F}_p\mathrm{Id}$ there exists a universal deformation ring; otherwise $R_{G_S}(\bar{\rho})$ will merely stand for a versal deformation ring.

## 3. *Strategy*

We use the notation $\rho$ for $\rho_G$ whenever no confusion shall arise.
The profinite group $G$ has a normal $p$-Sylow subgroup $P_G$ of finite index and finite type. After lemma 2.4 [Bo1] there exists exactly one semi-direct product $G = A \ltimes P_G$ with given action of $A$ on $\bar{P}_G$, where $\bar{P}_G = P_G/P_G^p[P_G, P_G]$. One then writes $\rho(G) = \rho(A) \ltimes \rho(P_G)$. To be more precise (proposition 2.3 [Bo1])

**Lemma 3.1** *If $V$ is an $\mathbf{F}_p[A]$-submodule of $\bar{P}_G$, then there exists an $A$-invariant subgroup $B$ of $P_G$ with $dim_{\mathbf{F}_p}V$ generators mapping onto $V$.*

Since $(\sharp A, p) = 1$ and $\bar{\rho} : A \to \mathrm{Gl}_2(\mathbf{F}_p)$ is injective, $\rho(A)$ is known by proposition 3.1 [Bo3]; the universal deformation ring of $\bar{\rho}_{|A}$ is $\mathbf{Z}_p$ and $\rho(A)$ is given by the canonical homomorphism $\mathbf{Z}_p \to R$.

Boston replaces the functor $\mathcal{F}$ by the functor $\mathcal{E} : \mathcal{C} \to \mathrm{Sets}$ defined by $\mathcal{E}(R) = \mathrm{Hom}_A(P_G, \Gamma_2(R))$ which is always representable. To be more specific:

- If $K = F$, $\mathcal{F} \cong \mathrm{Hom}_A(P_G, \Gamma_2(R))$.
- Otherwise $\rho(P_G) \not\subset \Gamma_2(R)$ but to compute the universal representation $\rho$ it suffices to know $\rho(A)$, $\mathrm{Hom}_A(\mathrm{Gal}(M/K), \Gamma_2(R))$, the images $\rho(\mathrm{Gal}(K/F))$ and $\rho(\mathrm{Gal}(M/K))$, as well as the action of $\rho(\mathrm{Gal}(K/F))$ on $\rho(\mathrm{Gal}(M/K))$.

**Remark 3.2** These considerations are also valid if we replace $G$ and $M$ by $G_S$ and $K_{S,p}$.

**Remark 3.3** The above considerations justify the method announced in the introduction: on the one hand we find the image of the generators of $P_G$ using the action of $A$; for

$x \in \mathrm{Ker}\bar{\rho}$, $\rho(x) = \begin{pmatrix} 1+Y_1 & Y_2 \\ Y_3 & 1+Y_4 \end{pmatrix}$, $Y_i \in \mathcal{M}_R$, $1 \le i \le 4$; this image introduces four variables. Expressing the action of $A$ allows to reduce this number of variables. We begin exactly as if $P_G$ were free. We choose a system of generators the images of which are particularly simple.

On the other hand we express the relations between selected generators. We need the relations of $P_G$. For this, we use that $\Gamma \cong P_G/X_\infty$ is free so we merely have to determine the relations of $X_\infty$ as a $\mathbf{Z}_p[[\Gamma]]$-module. These relations will yield relations between the variables in a straightforward fashion.

## 4. *Preliminary computations*

Following Boston's approach, we begin by expressing the action of $\rho(A)$ on $\rho(G)$. In particular, subsection 4.2 shows that $A$ possesses a uniquely determined complex conjugation $c$. This complex conjugation plays an important rôle because it allows us to reduce the number of variables. We then express the action of the other elements of $A$, if there are any. When $\rho(A)$ acts on a commutative group, its action is expressed in terms of basic linear algebra. We therefore find it convenient to work with $\bar{P}_G$.

4.1. *Image of the residual representation.* To express the action of $A$, information on the image of the residual representation $\bar{\rho}$ is needed. In the following lemma we show that $\mathrm{Gal}(K/F) \cong \bar{\rho}(P_G)$ and $A$ are abelian.

**Lemma 4.1.1** *One has* $\bar{\rho}(P_G) \subset \begin{pmatrix} 1 & * \\ 0 & 1 \end{pmatrix}$ *and $A$ is isomorphic to a subgroup of invertible diagonal matrices.*

Proof: The group $P_G$ is a pro-$p$-group, thus $\bar{\rho}(P_G)$ is a $p$-group which acts on the $p$-group $\mathrm{M}_2(\mathbf{F}_p)$, hence $\bar{\rho}(P_G) \subset \begin{pmatrix} 1 & * \\ 0 & * \end{pmatrix}$ (for a valuable basis). If $\begin{pmatrix} 1 & b \\ 0 & b' \end{pmatrix} \in \bar{\rho}(P_G)$, then there exists $s \in \mathbf{N}$ such that $b'^{p^s} = 1$ since $\mathrm{Gal}(K/F)$ is a $p$-group; yet $b' \in \mathbf{F}_p$, hence $b' = 1$ and

$$\bar{\rho}(P_G) \subset \begin{pmatrix} 1 & * \\ 0 & 1 \end{pmatrix}$$

One has $A \cong \mathrm{Im}\bar{\rho}/\bar{\rho}(P_G)$, hence $A$ is isomorphic to a subgroup of the group of invertible diagonal matrices. ∎

**Remark 4.1.2** As $\bar{\rho}$ is odd and $A$ is abelian, $A$ possesses a uniquely determined complex conjugation.

4.2. *Action of complex conjugation on $P_G$.* We can assume $\rho(c) = \begin{pmatrix} 1 & 0 \\ 0 & -1 \end{pmatrix}$ because we want $\rho$ up to strict equivalence.

**Lemma 4.2.1** *Let $x \in P_G$, $\rho(x) = \begin{pmatrix} U_1 & U_2 \\ U_3 & U_4 \end{pmatrix} \in \mathrm{Gl}_2(R)$, then*

- *If $x$ is invariant under complex conjugation then*

$$\rho(x) = \begin{pmatrix} 1+Y & 0 \\ 0 & 1+Y' \end{pmatrix} \quad \text{with } Y, Y' \in \mathcal{M}_R$$

- *If $x$ is inverted by complex conjugation then*

$$\rho(x) = \begin{pmatrix} (1+UU')^{1/2} & U \\ U' & (1+UU')^{1/2} \end{pmatrix}, \quad \text{with } U \in R, U' \in \mathcal{M}_R$$

Proof: It suffices to compare $\rho(x)$ or $\rho(x)^{-1}$ to

$$\rho(c \cdot x) = \rho(c)\rho(x)\rho(c)^{-1} = \begin{pmatrix} U_1 & -U_2 \\ -U_3 & U_4 \end{pmatrix}$$

and to use lemma 4.1.1 to find the values of the introduced coefficients. ∎

4.3. *Action of A.* If $x \in P_G$, $\bar{x}$ denotes the image of $x$ in $\bar{P}_G$. Since $A$ is abelian, the group $\bar{P}_G$ decomposes into $A$-invariant subgroups $< \bar{x} >$ of dimension 1. After lemma 3.1 there exists $x \in P_G$ that maps onto $\bar{x}$ and $\chi$, a character of $A$ such that

$$x \in P_{G,\chi} := \{u \in P_G : a \cdot u = u^{\chi(a)}, \forall a \in A\}$$

In the Borel case the action of $\rho(A)$ on $\rho(X)$ is known explicitly.
Recall that $\rho(A)$ is abelian; since it contains $\begin{pmatrix} 1 & 0 \\ 0 & -1 \end{pmatrix}$, one has
$\rho_{|A} = \begin{pmatrix} \chi_1 & 0 \\ 0 & \chi_2 \end{pmatrix}$ and:

**Proposition 4.3.1** *Let $x \in P_{G,\chi}$ with $\chi$ an odd character; then the action of $A$ on $x$ imposes:*

$$\rho(x) = \begin{pmatrix} 1 & U \\ 0 & 1 \end{pmatrix} \quad \text{if } \chi = \chi_1\chi_2^{-1}$$

$$\rho(x) = \begin{pmatrix} 1 & 0 \\ U' & 1 \end{pmatrix} \quad \text{if } \chi(a) = \chi_2\chi_1^{-1}$$

$$\rho(x) = \mathrm{Id} \quad \textit{otherwise.}$$

PROOF: Let $a \in A$; we denote by

$$\rho(a) = \begin{pmatrix} a_1 & 0 \\ 0 & a_4 \end{pmatrix}, \quad a_i \in R^*.$$

Recall that $\bar{\rho}(x) = \begin{pmatrix} 1 & b \\ 0 & 1 \end{pmatrix}$ with $b \in \mathbf{F}_q$. One writes $\rho(a \cdot x) = \rho(x)^{\chi(a)}$ and one uses the invariance of the trace by conjugation in order to obtain that if $\chi(a) \neq \pm 1$ then $UU' = 0$ (for details see theorem 3.8 [Bo3]). Then the action of $a$ on $x$ imposes one condition out of the four following ones:

- Prime-to-adjoint condition $U = U' = 0$

  if $\chi(a) = \pm 1$ and $a_1 \neq \pm a_4$,

if $\chi(a) \neq \pm 1$ and $a_1 \neq \chi(a)^{\pm 1} a_4$.

- Condition-free

  if $\chi(a) = \pm 1$ and $a_1 = \pm a_4$.

- $U' = 0$

  if $\chi(a) \neq \pm 1$, and $a_1 = \chi(a) a_4$.

- $U = 0$

  if $\chi(a) \neq \pm 1$ and $a_4 = \chi(a) a_1$.

Recalling that $\bar{\rho}_{|A} \neq \begin{pmatrix} \chi_1 & 0 \\ 0 & \pm\chi_1 \end{pmatrix}$, there must exist $a \in A$ such that $\rho(a) = \begin{pmatrix} a_1 & 0 \\ 0 & a_4 \end{pmatrix}$ with $a_1 \neq \pm a_4$, which shows the lemma in all cases. ∎

One notices that the conditions imposed by the action of $A$ reduce the number of deformation variables, but introduce no relations between these variables.

**Remark 4.3.2** In subsection 3.B [Bö1] Böckle solved the Borel case for an even representation ($\bar{\rho}(c) = \pm \mathrm{Id}$). If $\bar{\rho}$ is odd, the deformations can introduce such images as

$$\begin{pmatrix} (1+UU')^{1/2} & U \\ U' & (1+UU')^{1/2} \end{pmatrix}$$

which render the computation inextricable. To avoid this ugly image, we have imposed

$$\bar{\rho}_{|A} \neq \begin{pmatrix} \chi_1 & 0 \\ 0 & \pm\chi_1 \end{pmatrix}$$

If $x \in P_{G,\chi}$ with $\chi$ even, we obtain in the same fashion

**Proposition 4.3.3** *Let $x \in P_{G,\chi}$ with $\chi$ even; then*

$$\rho(x) = \begin{pmatrix} 1+Y & 0 \\ 0 & 1+Y' \end{pmatrix} \quad Y, Y' \in \mathcal{M}_R$$

*Moreover if* $\chi(A) \neq 1$, $\rho(x) = \mathrm{Id}$.

If $\Gamma = \Gamma_{cyc} =< \gamma >$ then $\bar{\gamma}$ (image of $\gamma$ in $\bar{P}_G$) is invariant by the action of $A$. Hence one can lift it to an element $s_d$ in $P_{G,triv}$.

**Lemma 4.3.4** *One has*

$$\rho(s_d) = \begin{pmatrix} 1+Y & 0 \\ 0 & 1+Y' \end{pmatrix} \text{ with } Y, Y' \in \mathcal{M}_R$$

PROOF: The element $s_d \in \mathrm{Ker}\ \bar{\rho}$ is invariant by the action of $A$; the expression of the action of the complex conjugation allows us to conclude. ∎

4.4. *Choice of representatives of the strict equivalence classes.* The chosen form above for the image of complex conjugation under $\rho$ does not completely fix a unique representative of a strict equivalence class. The representation $\rho$ can be replaced by its conjugate $\rho' = M\rho M^{-1}$ for some $M \equiv \mathrm{Id} \bmod \mathcal{M}_R$ without changing the strict equivalence class of $\rho$ provided $M$ commutes with $\rho(c)$. The matrix $M$ can be chosen in the form (see lemma 4.2.1)

$$M = \begin{pmatrix} 1+Z_1 & 0 \\ 0 & 1+Z_2 \end{pmatrix} \quad Z_1, Z_2 \in \mathcal{M}_R$$

This remark allows us to establish

**Lemma 4.4.1** *If* $\mathrm{Im}\bar{\rho}$ *is not diagonal, there exists* $x_n \in P_G$ *such that we can impose*

$$\rho(x_n) = \begin{pmatrix} 1 & 1 \\ 0 & 1 \end{pmatrix}$$

*Moreover if* $x \in P_{G,\chi}$ *commutes with* $x_n$ *and*

-*if* $\chi = \chi_1\chi_2^{-1}$ *then* $\rho(x) = \begin{pmatrix} 1 & U \\ 0 & 1 \end{pmatrix}$

-*if* $\chi = 1$ *then* $\rho(x) = \begin{pmatrix} 1+Y & 0 \\ 0 & 1+Y \end{pmatrix}$

-$\rho(x) = \mathrm{Id}$ *otherwise.*

Proof: If $\mathrm{Im}\bar{\rho}$ is not diagonal, there exists $x_n$ such that $\bar{\rho}(x_n) = \begin{pmatrix} 1 & 1 \\ 0 & 1 \end{pmatrix}$. Hence $x_n \in P_G$ is inverted by $c$ and using proposition 4.3.1,

$$\rho(x_n) = \begin{pmatrix} 1 & U_n \\ 0 & 1 \end{pmatrix} \text{ where } U_n \equiv 1 \bmod \mathcal{M}_R$$

Take $M$ as before. One has

$$M\rho(x_n)M^{-1} = \begin{pmatrix} 1 & \frac{1+Z_1}{1+Z_2}U_n \\ 0 & 1 \end{pmatrix}$$

As $U_n$ is invertible one can impose $U_n = 1$ by fixing the representative of the strict equivalence class of $\rho$.

The conditions of commutativity with $x_n$ give the results announced in the lemma. ∎

5. *Images of the generators*

5.1. *Minimal system of generators.* Let $d = \dim_{\mathbf{F}_p} H^1(G, \mathbf{F}_p) = \dim_{\mathbf{F}_p} H^1(P_G, \mathbf{F}_p)$ and $k = \dim_{\mathbf{F}_p} H^1(\Gamma, \mathbf{F}_p)$.

**Lemma 5.1.1** *A minimal system of ($\Lambda = \mathbf{Z}_p[[\Gamma]]$)-generators of the $\Lambda$-module $X_\infty$ contains $n = d - k$ elements.*

Proof: We write the inflation-restriction sequence of

$$0 \to X_\infty \to P_G \to \Gamma \to 0$$

using $\Gamma$ is free. The obtained exact sequence allows to conclude

$$0 \to H^1(\Gamma, \mathbf{F}_p) \to H^1(P_G, \mathbf{F}_p) \to H^1(X_\infty, \mathbf{F}_p)^\Gamma \to 0$$

∎

In the next paragraph we fix some notations.

If $X$ is a $p$-group and $x \in X$, $\bar{X}$ denotes the quotient of $X$ by its Frattini subgroup and $\bar{x}$ the image of $x$ in $\bar{X}$.

As $\bar{\Gamma}$ and $\overline{(X_\infty)_\Gamma}$ are invariant by $A$, we can decompose $\bar{P}_G$ as an $\mathbf{F}_p[A]$-module into:

$$\bar{P}_G = \overline{(X_\infty)_\Gamma} \oplus \bar{\Gamma}$$

To be more specific $\bar{\Gamma}$ is a $\mathbf{F}_p$-vector subspace of $\bar{P}_G$ with a $\mathbf{F}_p$-basis of eigenvectors for the action of $A$ $(\bar{s}_{n+1}, \ldots, \bar{s}_d)$ which admits as a complementary vector space $\overline{(X_\infty)_\Gamma}$. Let $(\bar{s}_1, \ldots \bar{s}_n)$ be a basis of $\overline{(X_\infty)_\Gamma}$ composed of eigenvectors for the action of $A$. Using lemma 3.1 we can lift these vectors to a minimal system $(s_1, \ldots, s_d)$ of generators of $P_G$.

We denote by $\tilde{\Gamma}$ the closed subgroup of $P_G$ generated by $(s_{n+1}, \ldots, s_d)$.

With these notations lemma 5.1.1 reads: since the quotient $P_G/X_\infty = \Gamma$ is pro-$p$-free, there exists a section $\Gamma \to P_G$ with image $\tilde{\Gamma}$ and we have

$$P_G = \tilde{\Gamma} \ltimes X_\infty$$

5.2. *Definition of the special generators.* By the previous subsection, there exists some minimal systems of generators of $\tilde{\Gamma}$, $P_G$ and a minimal system of $\Lambda$-generators of $X_\infty$, for which we control the action of $A$. Thus we know the image under $\rho$ of such generators. In this subsection, we rearrange these sytems according to the images by $\rho$.

**Definition 5.2.1** If $\mathrm{Im}\bar{\rho}$ is not diagonal, for the subgroup $X_\infty$ of $P_G$, there exists a minimal system of generators $(s_1, \ldots, s_n)$ and some integers $u_{X_\infty}, v_{X_\infty} \in \mathbf{N}$ such that:

- $\rho(s_i) = \begin{pmatrix} 1+Y & 0 \\ 0 & 1+Y \end{pmatrix}$, $Y \in \mathcal{M}_R$ , $1 \leq i \leq u_{X_\infty}$
- $\rho(s_i) = \begin{pmatrix} 1 & Y \\ 0 & 1 \end{pmatrix}$, $Y \in \mathcal{M}_R$, $u_{X_\infty} + 1 \leq i \leq v_{X_\infty}$
- $\rho(s_i) = \mathrm{Id}$, $v_{X_\infty} + 1 \leq i < n$
- $\rho(s_n) = \begin{pmatrix} 1 & 1 \\ 0 & 1 \end{pmatrix}$

These generators are termed *special generators*.

**Remark 5.2.2** If $(s_1 \ldots, s_n)$ is a minimal system of special generators then $s_i \in P_{G,\chi'_i}$ where $\chi'_i = \mathrm{Id}$ if $1 \leq i \leq u_{X_\infty}$, $\chi'_i = \chi_1\chi_2^{-1}$ if $u_{X_\infty}+1 \leq i \leq v_{X_\infty}$ or $i = n$ (see propositions 4.3.1 and 4.3.3).

**Definition 5.2.3** If $\mathrm{Im}\bar{\rho}$ is not diagonal, let $(s_1, \ldots, s_d)$ be a minimal system of generators of $P_G$ such that:

-the system $(s_1, \ldots, s_n)$ is a minimal system of special $\Lambda$-generators of $X_\infty$ and

-the system $(s_{n+1}, \ldots, s_d)$ is a minimal system of $\tilde{\Gamma}$ and $u_\Gamma, v_\Gamma, w_\gamma \in \mathbf{N}$ such that

- $\rho(s_{n+i}) = \begin{pmatrix} 1+Y & 0 \\ 0 & 1+Y' \end{pmatrix}$, $Y \neq Y' \in \mathcal{M}_R$ , $1 \leq i \leq u_\Gamma$
- $\rho(s_{n+i}) = \begin{pmatrix} 1 & Y \\ 0 & 1 \end{pmatrix}$, $Y \in \mathcal{M}_R$, $u_\Gamma + 1 \leq i \leq v_\Gamma$
- $\rho(s_{n+i}) = \begin{pmatrix} 1 & 0 \\ Y & 1 \end{pmatrix}$, $Y \in \mathcal{M}_R$, $v_\Gamma + 1 \leq i \leq w_\Gamma$
- $\rho(s_{n+i}) = \mathrm{Id}$, $w_\Gamma + 1 \leq i < k$

The generators $(s_{n+1}, \ldots, s_d)$ of $\tilde{\Gamma}$ are termed *special generators* of $\tilde{\Gamma}$.

The generators $(s_1, \ldots, s_d)$ of $P_G$ are termed *special generators* of $P_G$.

At last, a minimal system of generators of $P_{G_S}$, the image of which in the quotient $P_G$ is a minimal system of special generators of $P_G$ is termed minimal system of *special generators* of $P_{G_S}$.

**Remark 5.2.4** Since $\mathrm{Gal}(F/\mathbf{Q})$ is abelian, Leopoldt's conjecture holds for $F$. It implies that $u_\Gamma \leq 1$.

We denote by $Y_i$ the variable introduced by $\rho(s_i)$, $1 \leq i \leq v_{X_\infty}$; we denote by $Y_{2i-1+v_{X_\infty}}$ and $Y_{2i+v_{X_\infty}}$ the variables introduced by

$$\rho(s_{n+i}) = \begin{pmatrix} 1 + Y_{2i-1+v_{X_\infty}} & 0 \\ 0 & 1 + Y_{2i+v_{X_\infty}} \end{pmatrix},\ 1 \leq i \leq u_\Gamma$$

and by $Y_{i+u_\Gamma+v_{X_\infty}}$ the variable introduced by $\rho(s_{i+n})$, $u_\Gamma + 1 \leq i \leq w_\Gamma$.

**Remark 5.2.5** If $\mathrm{Im}\bar{\rho}$ is diagonal, we can define the special generators of $X_\infty$ as for $\tilde{\Gamma}$ and assume that $w_{X_\infty} = v_{X_\infty}$ in order to recover an situation analogous to the previous one.

## 6. *Determination of the universal deformation ring*

The study of sections 3-4-5 allows us to choose a minimal system of special generators of $P_G$, the images of which are simple. Since $P_G$ is not always free, we have to express the image by $\rho$ of the relations between selected generators in order to obtain the ideal of relations $I$ of $R_G(\bar{\rho})$. For this, we use the decomposition $P_G \cong \tilde{\Gamma} \ltimes X_\infty$.

6.1. *Action of* $\Gamma$ *on* $X_\infty$. *Choice of* $\Gamma$. Let $\Gamma$ be pro-$p$-free on $k$ generators. Let $\pi$ denote the projection $P_G \to \Gamma$ and $(\pi(s_{n+1}), \ldots, \pi(s_{n+k}))$ the image in $\Gamma$ of a minimal system of special generators of $\tilde{\Gamma}$. It is known that the correspondance that maps $\pi(s_{n+i})$ onto $1 + T_i$ induces an isomorphism of $\mathbf{Z}_p$-algebra $\Lambda = \mathbf{Z}_p[[\Gamma]] \cong \mathbf{Z}_p[[T_1, \ldots, T_k]]_{nc}$: this is the so called *Magnus algebra* which is non commutative except for $k = 1$, in which case it coincides with the usual Iwasawa algebra.
We can express the image under $\rho$ of the action of $\Lambda$ on $X_\infty$:

$$\rho(T_i x) = \rho((s_i - 1) \cdot x) = \rho(s_i)\rho(x)\rho(s_i)^{-1}\rho(x)^{-1}, \;\; n+1 \le i \le d, \; x \in X_\infty$$

To be more specific:

**Lemma 6.1.1** *Let $s$ (resp. $x$) be a special generator of $\tilde{\Gamma}$ (resp. $X_\infty$);*

*(i) if* $\rho(s) = \begin{pmatrix} 1+Y & 0 \\ 0 & 1+Y' \end{pmatrix}$ *and if* $\rho(x) = \begin{pmatrix} 1 & U \\ 0 & 1 \end{pmatrix}$ *then*

$$\rho((s-1)\cdot x) = \begin{pmatrix} 1 & (\frac{1+Y}{1+Y'} - 1)U \\ 0 & 1 \end{pmatrix}$$

*(ii) if* $\rho(s) = \begin{pmatrix} 1 & 0 \\ Y & 1 \end{pmatrix}$ *and if* $\rho(x) = \begin{pmatrix} 1 & U \\ 0 & 1 \end{pmatrix}$ *then*

$$\rho((s-1)\cdot x) = \begin{pmatrix} 1 - YU & YU^2 \\ -Y^2U & Y^2U^2 + YU + 1 \end{pmatrix}$$

*(iii)* $\rho((s-1)\cdot x) = \mathrm{Id}$ *otherwise.*

This lemma shows that the action of $\Lambda$ on $X_\infty$ is very simple provided that $w_\Gamma = v_\Gamma$ i.e for all special generators $s$ of $\tilde{\Gamma}$

$$\rho(s) \neq \begin{pmatrix} 1 & 0 \\ Y & 1 \end{pmatrix}$$

If $\Gamma = \Gamma_{cyc}$ then $w_{\Gamma_{cyc}} = v_{\Gamma_{cyc}}$.

We first assume that $\Gamma = \Gamma_{max}$ and $(s_{n+1}, \ldots, s_d)$ is a system of special generators of $\tilde{\Gamma}_{max}$. We define $v_{\Gamma_{max}}, w_{\Gamma_{max}}$ as in subsection 5.2. Let $\tilde{\Gamma}'$ be a subgroup of the subgroup of $\tilde{\Gamma}_{max}$ isomorphic to the quotient of $\Gamma_{max}$ by the normal subgroup generated by $(s_{n+v_{\Gamma_{max}}+1}, \ldots, s_{n+w_{\Gamma_{max}}})$.
It suffices also to replace $G$ by $G'$ defined through $\pi(\tilde{\Gamma}')$, where $G'$ is defined by $\pi(\tilde{\Gamma}')$ in the same way as $G$ is by $\pi(\tilde{\Gamma})$.

6.2. *Expression of the ideal of relations.* We first need to express the relations of $R_G(\bar{\rho})$ in terms of the relations of the ($\Lambda = \mathbf{Z}_p[[\Gamma]]$)-module $X_\infty$. We shall propose here a generalization of Boston's method ([Bo2] subsection 9.3).

Let $(f^i_j)_{1\le i\le n, 1\le j\le \ell}$ be the matrix of relations for a chosen system of special generators of the $\Lambda$-module $X_\infty$, where $\ell$ denotes the number of relations of $X_\infty$ for the above chosen system of special generators. To be more precise, the matrix of relations defined the endomorphism $\phi$ such that

$$\Lambda^\ell \xrightarrow{\phi} \Lambda^n \to X_\infty \to 0$$

**Proposition 6.2.1** *Let* $(s_1, \ldots, s_d)$ *a minimal system of special generators of* $P_G$ *(see subsection* 5.2. *for notations* $w_\Gamma, v_\Gamma, u_\Gamma, u_{X_\infty}, v_{X_\infty}$*). If* $w_\Gamma = v_\Gamma$ *and* $u_\Gamma = 1$*, (i.e* $\Gamma$ *contains the cyclotomic* $\mathbf{Z}_p$*-extension), then* $R_G(\bar{\rho}) = \mathbf{Z}_p[[Y_1, \ldots, Y_{d'}]]/I$*,* $I$ *being the ideal of relations generated by*

$$\prod_{i=1}^{u_{X_\infty}} (1+Y_i)^{f_i^j(0)} - 1,\ 1 \le j \le \ell$$

*and for* $1 \le j \le \ell$,

$$f_n^j\Big(\frac{1+Y_{v_{X_\infty}+1}}{1+Y_{v_{X_\infty}+2}} - 1, 0, \dots, 0\Big) +$$

$$\sum_{i=1+u_{X_\infty}}^{v_{X_\infty}} f_i^j\Big(\frac{1+Y_{v_{X_\infty}+1}}{1+Y_{v_{X_\infty}+2}} - 1, 0, \dots, 0\Big) Y_i$$

*where* $(\sum_{i=1}^n f_i^j(T_1, \dots, T_k)s_i)_{1 \le j \le \ell}$ *is a system of relations of the* $\Lambda$*-module* $X_\infty$*.*

PROOF: As in [Bo2] subsection 9.3 apply $\rho$ to the relations of the $\Lambda$-module $X_\infty$, using lemma 6.1.1 (iii) and (i). ∎

6.3. *Computation of the relations of the* $\Lambda$*-module* $X_\infty$*.* The relations of the $\Lambda$-module $X_\infty$ are in general not known. We shall give a way to derive them from the relations of $P_{G_S}$. To this end we shall proceed as follows: first the Fox derivative allows one to describe an intermediate $\Lambda$-module $Y_\infty$ with the help of generators and relations; second the diagram of lemma 6.3.1. yields a description of $X_\infty$ from that of $Y_\infty$.
We want to emphasize that the ideal of relations of $R_G(\bar{\rho})$ results directly from the relations of $P_{G_S}$ along a systematic line.

In [Ng1] (proposition 1.7) Nguyen Quang Do proposes a description of the $\Lambda$-module $Y_\infty = H_0(\mathrm{Gal}(K_{S,p}/F_\infty), I(K_{S,p}/F))$ based upon the exact sequence:

$$0 \to \Delta(F_\infty) \to \mathbf{Z}_p[[\mathrm{Gal}(F_\infty/F)]]^r \xrightarrow{\phi} \mathbf{Z}_p[[\mathrm{Gal}(F_\infty/F)]]^d \to Y_\infty \to 0$$

where $\phi$ denotes the composition of the projection of $\Gamma$ by the Fox derivative of the relations of $\mathrm{Gal}(M/F)$, and where $d = \dim_{\mathbf{F}_p} H^1(G_S, \mathbf{F}_p)$, and $r = \dim_{\mathbf{F}_p} H^2(G_S, \mathbf{F}_p)$. One also denotes by $I(K_{S,p}/F)$ the augmentation ideal of $\mathrm{Gal}(K_{S,p}/F)$.
Moreover one has $\Delta(F_\infty) = H^2(\mathrm{Gal}(K_{S,p}/F_\infty), \mathbf{Q}_p/\mathbf{Z}_p)^* = 0$ since, (see theorem 2.2 [Ng1]), Leopoldt's conjecture holds in view of $\mathrm{Gal}(F/\mathbf{Q})$ being abelian. The matrix of $\phi$ has rank $r - d$ and describes the relations of $Y_\infty$. In our case these relations give the relations of $X_\infty$ using the following lemma (see [Ja] lemma 4.3):

**Lemma 6.3.1** *In the commutative diagram below the two horizontal sequences and the two vertical sequences are exact:*

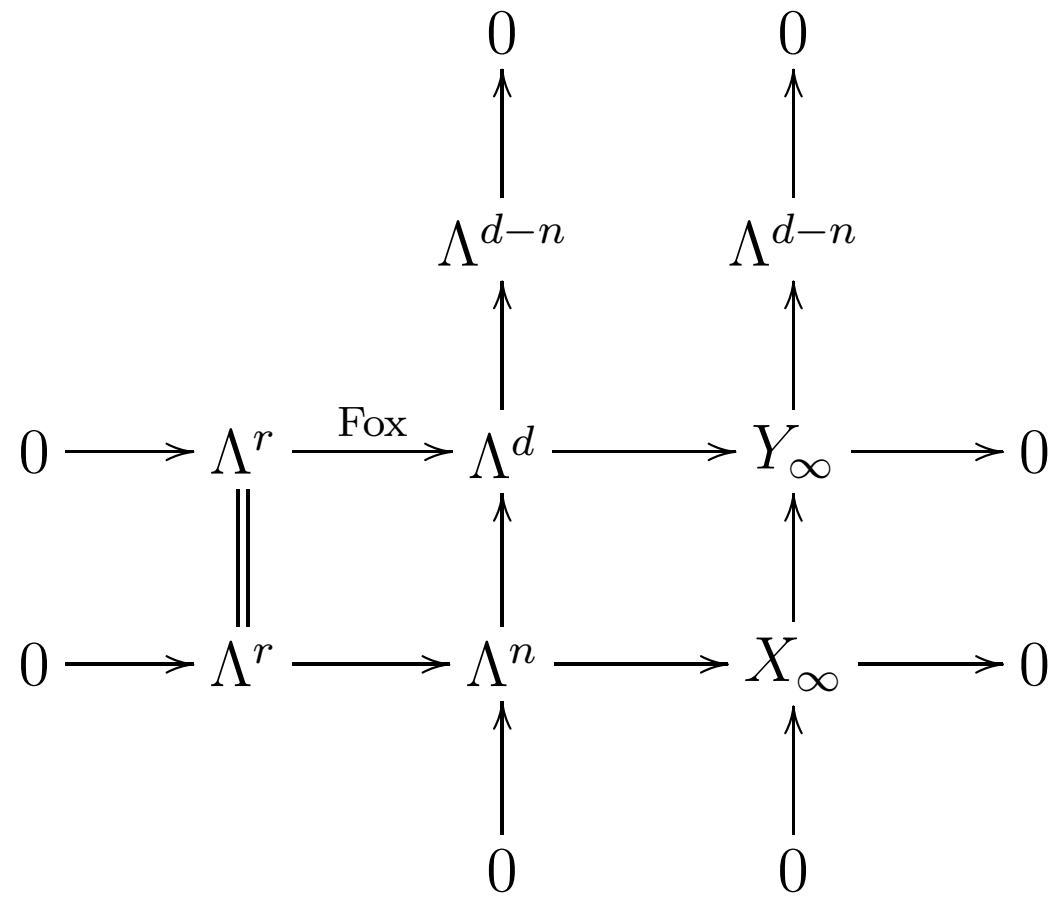

Proof: In the diagram

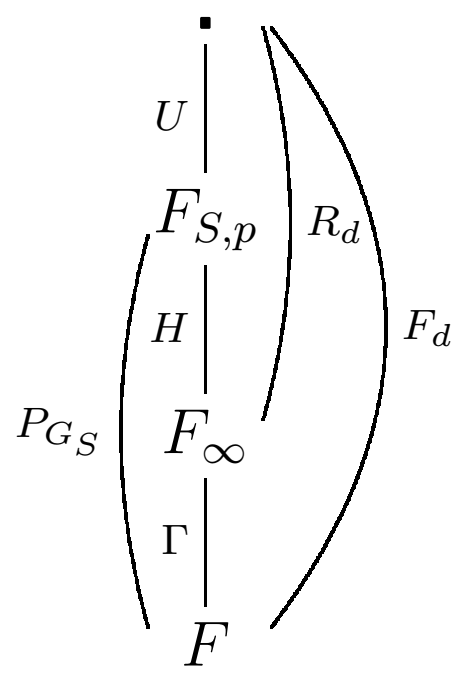

$F_d$ denotes a free presentation of rank $d$ of $P_{G_S}$; we assume $F_d$ is too a free presentation of $\Gamma$ and $R_d$ is defined by

$$0 \to R_d \to F_d \to \Gamma \to 0$$

The augmentation ideal of $F_d$ is defined by the exact sequence:

$$0 \to I(F_d) \to \mathbf{Z}_p[[F_d]] \to \mathbf{Z}_p \to 0$$

One applies to this exact sequence the functor $H_i(R_d, \centerdot)$. The freeness of $\mathbf{Z}_p[[F_d]]$ implies the Lyndon exact sequence ([Ng1] proposition 1.1):

$$0 \to H_1(R_d, \mathbf{Z}_p) \to H_0(R_d, I(F_d)) \to H_0(R_d, \mathbf{Z}_p[[F_d]]) \to H_0(R_d, \mathbf{Z}_p)$$

We deduce:

$$0 \longrightarrow R_d^{ab} \longrightarrow \Lambda^d \longrightarrow I(\Lambda) \longrightarrow 0$$

where $R_d^{ab} = R_d/[R_d, R_d]$. Similarly one applies the functor $H_i(H, \centerdot)$ to the exact sequence which defines the augmentation ideal of $P_{G_S}$, so that ([Ng1] proposition 1.7)

$$0 \to H_1(H, \mathbf{Z}_p) \to H_0(H, I(P_{G_S})) \to H_0(H, \mathbf{Z}_p[[P_{G_S}]]) \to H_0(H, \mathbf{Z}_p) \to 0$$

Hence,

$$0 \longrightarrow X_\infty \longrightarrow Y_\infty \longrightarrow I(\Lambda) \longrightarrow 0$$

The following commutative diagram

$$\begin{array}{ccc} R_d & \longrightarrow & H \\ \circlearrowright & & \circlearrowright \\ F_d & \longrightarrow & P_{G_S} \end{array}$$

yields functorial morphisms

$$H_0(R_d, I(F_d)) \to H_0 = (H, I(P_{G_S})) \text{ and } H_0(R_d, \mathbf{Z}_p[[F_d]]) \to H_0(H, \mathbf{Z}_p[[P_{G_S}]])$$

This proves the commutativity of the following diagram

$$\begin{array}{ccccc} 0 & & 0 & & \\ \uparrow & & \uparrow & & \\ I(\Lambda) & \longrightarrow & I(\Lambda) & & \\ \uparrow & & \uparrow & & \\ \Lambda^d & \longrightarrow & Y_\infty & \longrightarrow & 0 \\ \uparrow & & \uparrow & & \\ R_d^{ab} & \longrightarrow & X_\infty & \longrightarrow & 0 \\ \uparrow & & \uparrow & & \\ 0 & & 0 & & \end{array}$$

We shall complete the horizontal sequences by identifying them with known exact sequences.
By the equalities

$$\Lambda^d = \mathbf{Z}_p[[\Gamma]]^d = H_0(R_d, I(F_d)) = H_0(\mathrm{Gal}(F_{S,p}/F_\infty), H_0(U, I(F_d)))$$

the arrows between these modules and $Y_\infty$ identify by functoriality. The corresponding line is completed by the Fox derivative: recall that the morphism 'Fox' is obtained thanks to the exact sequences:

$$0 \to H_0(H, U^{ab}) \to H_0(H, H_0(U, I(F_d))) \cong \Lambda^d \to Y_\infty \to 0$$

The same reasoning applies to the horizontal line which then identifies to

$$0 \longrightarrow H_1(U, \mathbf{Z}_p) \longrightarrow R_d^{ab} \longrightarrow X_\infty \longrightarrow 0$$

where $R_d^{ab} = R_d/[R_d, R_d]$.

Recall $k = \dim_{\mathbf{F}_p} H^1(\Gamma, \mathbf{F}_p)$ we can use the corollary 1.2 of [Ng1], to show that $R_d^{ab} \cong R_k^{ab} \times \Lambda^{d-k}$, where $R_k^{ab} = R_k/[R_k, R_k]$ and $R_k$ is defined by a presentation of $\Gamma$:

$$0 \to R_k \to F_k \to \Gamma \to 0$$

Since $\Gamma$ is free, one knows ([Br] lemma 5.1) that $I(\Lambda)$ is a free $\mathbf{Z}_p[[\Lambda]]$-module of rank $k$. Hence $R_k^{ab} = 0$, which yields the bottom horizontal exact sequence:

$$H_0(H, U^{ab}) \cong \Lambda^r \text{ and } H_1(U, \mathbf{Z}_p) \cong \Lambda^r$$

and $I(\Lambda) \cong \Lambda^k = \Lambda^{d-n}$. ∎

6.4. *General presentation by generators and relations.* We now apply our results using a presentation of $P_{G_S}$ by generators and relations. The pro-$p$-group $P_{G_S}$ has a minimal presentation by the exact sequence:

$$0 \to R \to F_d \to P_{G_S} \to 0$$

where $F_d$ is a free pro-$p$-group of rank $d$ and $R$ is the pro-$p$-group of the relations of $P_{G_S}$.
Let $\{s_i,\ 1 \leq i \leq d\}$ be a minimal system of generators of $P_{G_S}$ and $\{r_j, 1 \leq j \leq r\}$ a minimal system of relations of $P_{G_S}$ .
In order to render the injective morphism $X_\infty \to Y_\infty$ explicit, it is necessary to choose judiciously the generators of $P_{G_S}$. To this end we notice that the elements of $X_\infty$ are characterized by their image being trivial in the quotient $\Gamma$.

As in subsection 6.1, $\pi$ denotes the projection $P_G \to \Gamma$ and if $x \in P_{G_S}$, we denote $\tilde{x}$ its projection on $P_G$.

Let $(s_1, \ldots, s_d)$ be a minimal system of special generators of $P_{G_S}$ such that: $\pi(\tilde{s}_i) = 1,\ 1 \leq i \leq n$ (recall that $(\tilde{s}_1, \ldots, \tilde{s}_n)$ is a minimal system of generators of $X_\infty$).
We can now compute the Fox matrix $M_\phi = \Big(\pi\Big(\dfrac{\partial r_j}{\partial s_i}\Big)\Big)_{1\leq j\leq l,\, 1\leq i\leq d}$, hence $Y_\infty$ as well. To obtain $X_\infty$ it suffices to omit the last $k$ lines of $M_\phi$, which gives the matrix $M_{\phi'}$ of $\phi'$ describing the exact sequence

$$0 \to \Lambda^r \overset{\phi'}{\to} \Lambda^n \to X_\infty \to 0$$

thus describing $X_\infty$ as a $\Lambda$-module.

**Proposition 6.4.1** *Let*

$$0 \to R \to F_d \to P_{G_S} \to 0$$

*be a minimal presentation of $P_{G_S}$, with $\{r_j,\ 1 \leq j \leq r\}$ a minimal system of relations and $\{s_i,\ 1 \leq i \leq d\}$ a minimal system of special generators such that $\pi(\tilde{s}_i) = 0,\ 1 \leq i \leq n$. Then the $\Lambda$-module $X_\infty$ is described by*

$$0 \to \Lambda^r \overset{\phi'}{\to} \Lambda^n \to X_\infty \to 0$$

*the matrix $M_{\phi'}$ of $\phi'$ being*

$$M_{\phi'} = \Big(\pi\Big(\frac{\partial r_j}{\partial s_i}\Big)\Big)_{1\leq i\leq n,\, 1\leq j\leq r}$$

**Remark 6.4.2** Recall that the injectivity of $\phi'$ is obtained thanks to the weak Leopoldt conjecture. In order to discuss Mazur's conjecture (see [Ma] subsections 1.6 and 1.10) we remark that this injectivity means that the relations of $X_\infty$ obtained in this way are independent.

## 7. *Explicit universal deformation ring*

7.1. *General computation of the Fox matrix.* We define by induction

$$P^{(0)} = P_{G_S} \text{ and } P^{(n+1)} = [P^{(n)}, P_{G_S}]$$

There exists a decomposition with coefficients in $\mathbf{Z}_p$, using the 'Hall collecting process':

$$r_j \bmod P^{(n+1)} \equiv$$

$$\prod_{i=1}^{d} s_i^{a_i^j} \prod_{1\le i<k\le d} [s_i, s_k]^{a_{i,k}^j} \cdots \prod_{1\le i_1<i_2\le d}^{i_3,\dots,i_n\in[1,d]} [\cdots[[s_{i_1}, s_{i_2}], s_{i_3}], \dots, s_{i_n}]^{a^j_{i_1,\cdots,i_n}}$$

**Remark 7.1.1** Let $p^l$ be the lcm of the orders $p^{l_i}$ or $\infty$ of the elements $s_i[P_{G_S}, P_{G_S}]$, $1 \le i \le d$, with the wild convention that $p^\infty = 0$.
The coefficients $a_i^j \bmod (p^l)^2$ and $a_{i,k}^j \bmod p^l$ are known ([Ko1] proposition 7.23). They are related to the transgression map $H^2(R, \mathbf{Z}/p^l\mathbf{Z}) \to H^1(P_{G_S}, \mathbf{Z}/p^l\mathbf{Z})$ and to the cup product of elements of $H^1(P_{G_S}, \mathbf{Z}/p^l\mathbf{Z})$.

The following two lemmas (which are straightforward) allow us to compute the Fox derivative of $r_j$.

**Lemma 7.1.2** *Let $u, v \in F_d$, $a \in \mathbf{N}^*$ then*

$$\frac{\partial(uv)}{\partial s} = \frac{\partial u}{\partial s} + u\frac{\partial v}{\partial s}, \quad \textit{hence} \quad \frac{\partial u^a}{\partial s} = \sum_{k=0}^{a-1} u^k \frac{\partial u}{\partial s}$$

**Lemma 7.1.3** *Let $u \in F_d$, $i, j \in [1, d]$, then*

$$\frac{\partial}{\partial s_i}[u, s_j] = \frac{\partial u}{\partial s_i} + us_j\frac{\partial u^{-1}}{\partial s_i} \ \textit{if } i \ne j$$

$$\frac{\partial}{\partial s_i}[u,s_i]=\frac{\partial u}{\partial s_i}+u+us_i\frac{\partial u^{-1}}{\partial s_i}-us_iu^{-1}s_i^{-1}$$

7.2. *The $\mathbf{Z}_p$-cyclotomic extension of $F$.* We make precise the computation of $R_G(\bar\rho_G)$ in the case $\Gamma=\Gamma_{cyc}=<\gamma>$ (recall $w_{\Gamma_{cyc}}=v_{\Gamma_{cyc}}=u_{\Gamma_{cyc}}=1$). We may assume that $\pi(s_d)=\gamma$, $\pi(s_i)=1,\ 1\le i\le d-1$. After lemma 7.1.3 the Fox matrix reads

$$M_\phi=\Big(\pi\Big(\frac{\partial r_j}{\partial s_i}\Big)\Big)=\Big(a_i^j+\sum_{k=1}^{\infty}b_{i,k}^j(\gamma-1)^k\Big)_{1\le i\le d,\ 1\le j\le r}$$

where $b_{i,k}^j$ denotes the coefficient of $[\cdots[[s_i,\gamma],\gamma],\dots,\gamma]\in P^{(k)}-P^{(k+1)}$ in the relation $r_j$.

Recall lemma 4.3.4 and the notations of subsection 5.2; for $\Gamma=\Gamma_{cyc}$ the variables introduced by $\tilde{s_d}$ are $Y_{d'-1},Y_{d'}$ in $R_G(\bar\rho_G)$ (see introduction: the number of variables of $R_G(\bar\rho)$ is $d'=\dim_{\mathbf{F}_p}H^2(G,Ad\bar\rho_G)$), hence:

$$\rho_G(\tilde{s}_d)=\begin{pmatrix}1+Y_{d'-1}&0\\0&1+Y_{d'}\end{pmatrix}\ \text{with}\ Y_{d'-1},Y_{d'}\in\mathcal{M}_R$$

We can also apply proposition 6.2.1 so as to obtain $R_G(\bar\rho_G)$:

**Theorem 7.2.1** *In the Borel case, if* $\mathrm{Im}\bar\rho$ *is not diagonal and if* $\Gamma=\Gamma_{cyc}$*, we obtain*

$$R_G(\bar\rho_G)=\mathbf{Z}_p[[Y_1,\dots,Y_{d'}]]/I$$

*where $I$ is the ideal of relations generated by*

$$\prod_{i=1}^{u_{X_\infty}}(1+Y_i)^{a_i^j}-1,\ 1\le j\le\ell$$

*and*

$$a_n^j+\sum_{k=1}^{\infty}b_{n,k}^j\Big(\frac{1+Y_{d'-1}}{1+Y_{d'}}-1\Big)^k+\sum_{i=1+u_{X_\infty}}^{v_{X_\infty}}\Big(a_i^j+\sum_{k=1}^{\infty}b_{i,k}^j\Big(\frac{1+Y_{d'-1}}{1+Y_{d'}}-1\Big)^k\Big)Y_i,\ 1\le j\le\ell$$

7.3. *Free pro-p-extension of rank k of F.* Let us assume $\Gamma$ has rank $k$. If $\Gamma_{cyc} =< \gamma >\subset \Gamma$ and $w_\Gamma = v_\Gamma$, the computation of $R_G(\bar{\rho}_G)$ is analogous to that of the previous case $\Gamma = \Gamma_{cyc}$. After lemma 6.1.1 the image by $\rho_G$ of the action of $\Lambda \cong \mathbf{Z}_p[[T_1, \ldots, T_k]]_{nc}$ is commutative, and even trivial if $i > 1$; that is, if $f \in \mathbf{Z}_p[[T_1, \ldots, T_k]]_{nc}$ is monomonial and $i > t_\Gamma$

$$\rho_G(T_i f x) = \text{ Id}$$

Hence only the commutators of the form $r = [\cdots[[s_i, \gamma], \ldots, \gamma]$ have a non trivial image by $r \mapsto \rho_G(\pi(\frac{\partial r}{\partial s_i})x)$.
The theorem 7.2.1 also holds in the latter case (note the coefficients $u_{X_\infty}, d'$ have changed) and gives a better approximation $R_G(\bar{\rho}_G)$ of $R_{G_S}(\bar{\rho})$.

**Remark 7.3.1** Our method can also give an approximation of the ideal $I$: if the relations of $P_{G_S}$ are known modulo $P^{n+1}$, theorem 7.2.1 allows us to control the approximation of $I$. Moreover, even if the computations are inextricable if $w_\Gamma \neq v_\Gamma$, we can also work modulo $P^{n+1}$ (for a given $n$ depending on our patience).

7.4. *Comparison between $R_{G_S}(\bar{\rho})$ and $R_G(\bar{\rho}_G)$.* The study carried out in section 3-4 also applies when $G$, $P_G$ are replaced by $G_S$, $P_{G_S}$. Since $\bar{P}_G = \bar{P}_{G_S}$, if $< \bar{x} >$ is an $A$-invariant subgroup, we can not only lift $\bar{x}$ to $\tilde{x} \in P_G$ but also to $x \in P_{G_S}$. The computations of 4.2 and 4.3 give the image of $\rho(x)$ as if $P_{G_S}$ were free.
Let us recall how to compute the universal ring in order to make precise the surjection $R_{G_S}(\bar{\rho}) \to R_G(\bar{\rho}_G)$.
After remark 3.1, if $\tilde{x} \in \text{Ker}\bar{\rho}$ then

$$\rho_G(\tilde{x}) = \begin{pmatrix} 1 + Y_1 & Y_2 \\ Y_3 & 1 + Y_4 \end{pmatrix}, \; Y_i \in \mathcal{M}_R, 1 \leq i \leq 4$$

This image a priori introduces four variables (depend upon $\tilde{x}$); the action of A allows to reduce the number of variables to

- 1 or 0 if $\tilde{x}$ is inverted by complex conjugation,
- 2 or 0 if $\tilde{x}$ is invariant by complex conjugation.

These computations apply to $\rho_G$ or $\rho$.

Taking $G$ instead of $G_S$ introduces additional conditions of commutativity which allow us to cancel the lower triangular images and to render scalar the diagonal images.

Then we express the image by $\rho_G$ (resp. $\rho$) of the relations of $P_G$ (resp. $P_{G_S}$). The universal deformation ring $R_G(\bar\rho_G)$ (resp. $R_{G_S}(\bar\rho)$) is the quotient of the ring of formal series in the introduced variables by the images of the special generators of $P_G$ (resp. $P_{G_S}$), by the image of the relations of $P_G$ (resp. $P_{G_S}$) by $\rho_G$ (resp. $\rho$).

The relations of $P_G$ are obtained from the relations of $P_{G_S}$ using the Fox derivative.

The pro-$p$-group $P_{G_S}$ has a minimal presentation by the exact sequence:

$$0 \to R \to F_d \to P_{G_S} \to 0$$

where $F_d$ is a free pro-$p$-group of rank $d$ and $R$ is the pro-$p$-group of the relations of $P_{G_S}$.

Let $\{s_i,\ 1 \leq i \leq d\}$ a minimal system of special generators of $P_{G_S}$.

Let $r$ be a relation of $P_{G_S}$; we can think of $r$ as a word in $F_d$. Let $\tilde r \in P_G$ be the element obtained from $r$ when we remplace the letters $s_i$ by $\tilde s_i$, $1 \leq i \leq d$.

We observe that

$$\rho_G(\tilde r) = \rho_G\Big(\sum_{i=1}^{n} \pi\Big(\frac{\partial r}{\partial s_i}\Big)\tilde s_i\Big)$$

Hence we obtain

**Theorem 7.4.1** *There exists a minimal system $\{s_i, 1 \leq i \leq d\}$ of (special) generators of $P_{G_S}$ and choice of the representatives of the strict equivalence classes $\rho$, $\rho_G$ such that the surjective morphism*

$$R_{G_S}(\bar\rho) \to R_G(\bar\rho_G)$$

*maps the variables introduced by $\rho(s_i)$ onto the variables introduced by $\rho_G(\tilde s_i)$.*

**Remark 7.4.2** The surjection $R_{G_S}(\bar\rho) \to R_G(\bar\rho_G)$ yields

$$\dim_{Krull} R_{G_S}(\bar\rho)/pR_{G_S}(\bar\rho) \geq \dim_{Krull} R_G(\bar\rho_G)/pR_G(\bar\rho_G)$$

This could allow us to discuss Mazur's question ([Ma] subsection 1.10).

7.5. *Examples.* **Cyclotomic fields**. The first examples of representations $\bar{\rho} : \mathrm{Gal}(\bar{\mathbf{Q}}/\mathbf{Q}) \to \mathrm{Gl}_2(\mathbf{F}_p)$ in the Borel case appear in the study of elliptic curves. These representations introduce the cyclotomic fields $F = \mathbf{Q}(\zeta_p)$, where $\zeta_p$ denotes a primitive $p$-root of unit.
We can find some examples of such representations associated to elliptic curves for $p = 5, 7$ in [Se] subsection 5.5. In these cases $R_{G_S}(\bar{\rho})$ is known because $p$ is regular hence $P_{G_S}$ is free.
In [Bo3] section 6, Boston gives an example of such a representation with $F = \mathbf{Q}(\zeta_p)$ for the irregular prime $p = 691$; this representation is associated to the unique normalized cusp form of weight 12.
We would like to correct an imprecision in [Bo2] proposition 9.2 and [Bo3] section 6, where a confusion seems to arise because 'Spiegelung' is overlooked (see below).
Let $\bar{\rho} : \mathrm{Gal}(\bar{\mathbf{Q}}/\mathbf{Q}) \to \mathrm{Gl}_2(\mathbf{F}_p)$ be a continuous odd representation unramified outside $S_{\mathbf{Q}} = \{p, \infty\}$ with $\mathrm{Im}\bar{\rho} = \begin{pmatrix} 1 & * \\ 0 & * \end{pmatrix}$ and $F = \mathbf{Q}(\zeta_p)$.
We denote by $a$ a generator of $A \cong \mathrm{Gal}(F/\mathbf{Q})$ and $\omega : A \to \mathbf{Z}_p^*$ the
Teichmüller lift of the cyclotomic character. There exists even $k' \in \mathbf{N}$ depending upon the representation such that

$$\bar{\rho}(a) = \begin{pmatrix} 1 & 0 \\ 0 & \omega^{k'-1}(a) \bmod p \end{pmatrix}$$

Assume that Vandiver's conjecture holds for $p$. Let $S_F$ be the set of places in $F$ over $S_{\mathbf{Q}}$, $\Gamma = \Gamma_{cyc}$ and $\tilde{\gamma}$ denotes a lift to $P_G$ of a generator of $\Gamma$.
Since we assume Vandiver's conjecture to hold, we know the structure of $X_\infty$ as a $\Lambda$-module.
For any $n \geq 1$, let $A_n$ be the $p$-class group of $\mathbf{Q}(\zeta_{p^n})$,

$$A_\infty = \varinjlim A_n, \ Z_\infty = \varprojlim A_n$$

Then:
-It is well known that Vandiver's conjecture implies $A_n^+ = 0$ for any $n \geq 1$, so that the natural maps $A_n \to A_\infty$ are injective. It follows that the $\Lambda$-torsion free part of $X_\infty$ is actually $\Lambda$-free, and hence we have an isomorphism of $\Lambda$-modules (and not only a pseudo-isomorphism)

$$X_\infty \cong \Lambda^{r_2} \oplus \mathrm{tor}_\Lambda X_\infty$$

where $r_2 = (p-1)/2$ is the number of complex places in $F$ and $\mathrm{tor}_\Lambda X_\infty$ denotes the $\Lambda$-torsion of $X_\infty$.

-The structure of $\mathrm{tor}_\Lambda X_\infty$ can be given by a 'mirror' ('Spiegelung') version of $Z_\infty$. More precisely, using the fact that each $\mathbf{Q}(\zeta_{p^n})$ admits only one place above $p$ and verifies Leopoldt's conjecture, we get (for details, see [Ng3] section 3):

$$\mathrm{tor}_\Lambda X_\infty \cong \lim_{\leftarrow} \mathrm{Hom}_{\Gamma_n}(Z_\infty, \mu_{p^\infty})$$

where $\Gamma_n = \mathrm{Gal}(\mathbf{Q}(\zeta_{p^\infty})/\mathbf{Q}(\zeta_{p^n}))$. Because $A_n^+ = 0$ for any $n \geq 1$, it follows that $X_\infty^+ = \mathrm{tor}_\Lambda X_\infty^+ = \mathrm{tor}_\Lambda X_\infty$, and that $X_\infty^-$ (which admits $Z_\infty$ as a quotient) is a $\Lambda$-free on $r_2 = (p-1)/2$ generators $(\tilde{x}_1, \ldots, \tilde{x}_{p-2})$ where $\tilde{x}_i \in P_{G,\omega^i}$, $1 \leq i \leq p-2$ (see [Wa] corollary 10.15).

Recall $A \cong \mathrm{Gal}(F/\mathbf{Q})$ is cyclic of order $p-1$ prime to $p$. Hence the subgroup of invariants by $A$ of $X_\infty$ is isomorphic to image of the norm $\nu$ of $X_\infty$:

$$X_\infty^A \cong \nu X_\infty$$

By class field theory, $\nu X_\infty$ is isomorphic to the analogue of $X_\infty$ for the cyclotomic $\mathbf{Z}_p$-extension $\mathbf{Q}_\infty/\mathbf{Q}$, hence is trivial. Recall that if $x \in X_{\infty,\omega^i}$ and $i$ even then $\rho(x) = \mathrm{Id}$ if $i \not\equiv 0 \bmod p$. Thus the image of $X_\infty^+$ by $\rho$ is trivial by using subsection 4.3 we get

**Proposition 7.5.1** *In the Borel case, if $\Gamma = \Gamma_{cyc}$, if $\mathrm{Im}\bar{\rho}$ not diagonal such that $F = \mathbf{Q}(\zeta_p)$ and $p < 1 + 4.10^6$ (so that Vandiver's conjecture holds), then*

$$R_G(\bar{\rho}_G) = \mathbf{Z}_p[[Y_1, Y_2]]$$

*and the universal representation is given by*

$$\rho_G(\tilde{a}) = \begin{pmatrix} 1 & 0 \\ 0 & \omega^{k'-1}(a) \end{pmatrix}, \ \rho_G(\tilde{\gamma}) = \begin{pmatrix} 1+Y_1 & 0 \\ 0 & 1+Y_2 \end{pmatrix}$$

$$\rho_G(\tilde{x}_{p-k'}) = \begin{pmatrix} 1 & 1 \\ 0 & 1 \end{pmatrix}$$

*and all other given generators of $P_G$ have a trivial image by $\rho_G$.*

Using theorem 7.4.1 we obtain

**Proposition 7.5.2** *In the Borel case, if $Im\bar{\rho}$ not diagonal such that $F = \mathbf{Q}(\zeta_p)$ and $p < 1 + 4.10^6$ (so that Vandiver's conjecture holds), then*

$$R_{G_S}(\bar{\rho}) = \mathbf{Z}_p[[Y_1, Y_2, Y_3]]/I$$

*The universal representation (with obvious notations) is given by*

$$\rho(a) = \begin{pmatrix} 1 & 0 \\ 0 & \omega^{k'-1}(a) \end{pmatrix}, \; \rho(\gamma) = \begin{pmatrix} 1+Y_1 & 0 \\ 0 & 1+Y_2 \end{pmatrix}$$

$$\rho(x_{p-k'}) = \begin{pmatrix} 1 & 1 \\ 0 & 1 \end{pmatrix}, \; \rho(x_{k'-1}) = \begin{pmatrix} 1 & 0 \\ Y_3 & 1 \end{pmatrix}$$

*and all other given generators of $P_{G_S}$ have a trivial image by $\rho$. Moreover if $r \in I$ we have $r \equiv 0$ mod $Y_3$, and $I = (0)$ if $p$ is regular.*

**Remark 7.5.3** For these representations, we have

$$\dim_{Krull} R_{G_S}(\bar{\rho})/pR_{G_S}(\bar{\rho}) \leq 3.$$

A reasoning analogous to that of [Ma] subsection 1.10 proves that the latter Krull dimension is at least 3 ([Bo1] remark 5.1). Thus this dimension is 3 and Mazur's question ([Ma] subsection 1.6) still holds in these cases.

**Increasing the ramification**. The problem of increasing the ramification has been studied by Boston [Bo4] theorem 1 and Böckle [Bö1] subsection 3.C (for even representations). Let $p$ be a regular odd prime, $F = \mathbf{Q}(\zeta_p)$ and $S_p$ be the set of places over $p$ and archimedean places.

In this case $G_{S_p} = \Gamma_{max}$ is uniquely determined by [Ya2] proposition 2.2 and is free with $(p+1)/2$ generators (see previous description). Remark that the quotient $\Gamma_{max}$ of $G_S$ does not depend on $S$ with $S_p \subset S$.

If we increase the ramification (i.e we consider $S$ instead of $S_p$), Neumann [Ne] corollary 5.3 gives a presentation of $G_S$:

**Lemma 7.5.4** *Let $(x_i,\ i \in I)$ be a minimal system of topological generators of $G_{S_p}$ and $(s_i,\ i \in I)$ be a system of elements of $G_S$ with*

$$s_i \ mod \ \mathrm{Gal}(K_S/K_{S_p}) = x_i,\ i \in I.$$

*Further let $t_q \in \mathrm{Gal}(K_S/K_{S_p})$ for $q \in S - S_p$ with $N(q) \equiv 1 \ mod \ p$, such that $t_q$ generates the inertia group of some prolongation of $q$.*

*Then the set*

$$\{s_i, t_q,\ i \in I,\ q \in S - S_p \ N(q) \equiv 1 \ mod \ p\}$$

*forms a set of generators of $G_S$ where the subset $\{t_q\}$ is free.*

Hence $\tilde{\Gamma}_{max}$ is generated by $(s_i, i \in \{0\} \cup \{1, 3, \ldots, (p-1)/2\})$ with $s_i \in P_{G_{S,\omega_i}}$.

Let $\bar{\rho} : G_S \to \mathrm{Gl}_2(F_p)$ in the Borel case, such that $F = \mathbf{Q}(\zeta_p)$ and there exists $q_0 \in S - S_p$ such that

$$\bar{\rho}(t_{q_0}) = \begin{pmatrix} 1 & 1 \\ 0 & 1 \end{pmatrix}$$

(see remark 2.1).

After the previous description, in order to have $w_\Gamma = 0$ (see subsection 6.1), we choose $\tilde{\Gamma} = \tilde{\Gamma}_{max}/ < s_{k'-1} >_{normal}$ (with $k'$ as in the previous paragraph). Hence we have the universal representation:

$$\rho_G(\tilde{s}_0) = \begin{pmatrix} 1+Y_1 & 0 \\ 0 & 1+Y_2 \end{pmatrix},\ \rho_G(\tilde{s}_{p-k'}) = \begin{pmatrix} 1 & Y_3 \\ 0 & 1 \end{pmatrix}$$

and $\rho_G(\tilde{s}_i) = \mathrm{Id}$ for the other $i \in I$;

$$\rho_G(\tilde{t}_{q_0}) = \begin{pmatrix} 1 & 1 \\ 0 & 1 \end{pmatrix}$$

$$\rho_G(\tilde{t}_q) = \begin{pmatrix} 1 & Y_q \\ 0 & 1 \end{pmatrix} \text{ or } \rho_G(\tilde{t}_q) = \begin{pmatrix} 1+Y_q & 0 \\ 0 & 1+Y_q \end{pmatrix}$$

where $q \in S - (S_p \cup \{q_0\})$ and some variables $Y_q$ can vanish depending of the action of $A$. In this particular situation we can easily apply the generalization of theorem 7.2.1 in subsection 7.3.

## 8. *Wingberg's presentation*

Recall the notations of the section 2: $\bar{\rho} : \mathrm{Gal}(\bar{\mathbf{Q}}/\mathbf{Q}) \to \mathrm{Gl}_2(\mathbf{F}_p)$ is an odd continuous representation in the Borel case unramified outside a finite set of primes $S_{\mathbf{Q}}$ of $\mathbf{Q}$, $K$ denotes the subfield of $\bar{\mathbf{Q}}$ fixed by $\mathrm{Ker}\bar{\rho}$, and $F$ is a subextension of $K$ such that $\mathrm{Gal}(K/F)$ is the Sylow $p$-subgroup of $\mathrm{Gal}(K/\mathbf{Q})$. The finite set $S_F$ of places of $F$ containing the places over $S_{\mathbf{Q}}$ in $F$, allows us to define $G_S$ and $P_{G_S}$; from here on $S$ denotes $S_F$.

8.1. *Using Wingberg's explicit presentation.* Our method is most effective whenever the relations of $P_{G_S}$ are known. Explicit results are known (see [Bo1]) in the case where $P_{G_S}$ is free. We shall now solve in a simple way the case where the relations of $P_{G_S}$ do not contain any threefold commutator. Then we can compute $R_{G_S}(\bar{\rho})$. A particular case of this corresponds to what we call a *Wingberg presentation*

Let $F_v$ denote the $v$-completion of $F$ for $v \in S$.

Let $\mu_p$ denote the set of $p$-roots of unity.

Let $\delta_v = \begin{cases} 1 & \text{if } \mu_p \subset F_v \\ 0 & \text{otherwise} \end{cases}$ and $\delta = \begin{cases} 1 & \text{if } \mu_p \subset F \\ 0 & \text{otherwise.} \end{cases}$

Let $S_0 \subset S$ be a maximal subset of finite primes such that:

$$\sum_{v \in S_0} \delta_v = \delta$$

Let

$$V_{S_0}^S = \{\alpha \in F^* \mid \alpha \in F_v^{*p} \text{ for } v \in S_0,\ \alpha \in U_v F_v^{*p},\ v \notin S\}/F^{*p}$$

where $U_v$ denotes the group of units in the ring of integers of $F_v$.

In [Wi] Wingberg proves the following proposition:

**Proposition 8.1.1** *The condition $V_{S_0}^S = 0$ is equivalent to $P_{G_S}$ being the free product of the decomposition groups $\mathcal{P}_v$, $v \in S - S_0$ and a free pro-p-group of rank $r_f = 1 + \sum_{v \in S_p \cap S_0} [F_v : \mathbf{Q}_p] - \sharp(S - S_0)$.*

**Remark 8.1.2** Leopoldt's conjecture for $F$ implies that the rank $k$ of $\Gamma_{max}$ verifies $1 \leq k \leq r_2 + 1$ where $r_2$ is the number of complex places in $F$. In the situation of proposition 8.1.1, Yamagishi [Ya1] has determined $k$ explicitely in terms of local datas; in particular $k \leq r_2 + 1$ (see also subsection 8.3).

**Remark 8.1.3** If $\mu_p \not\subset F$ then the decomposition groups $\mathcal{P}_v$, $v \in S - S_0$ are free; hence $P_{G_S}$ is free. The latter case is well known, so we now turn to $\mu_p \subset F$; in this case $S_0$ contains only one element.

If $p \neq 2$ then places that ramify in a pro-$p$-extension of a number field are primes $q$ of the form $N(q) \equiv 1 \bmod p$ or $q$ divides $p$. So we assume that $S_{\mathbf{Q}}$ and $S_F$ contain only places of this type.

Let $S'$ denote a set of $r_f$ places different from the elements of $S$. Hence $P_{G_S}$ is described

- by *free generators* $s_v$, $v \in S'$,
- by *tame generators* $s_v, t_v$ $v \in S - S_0$ such that $v \equiv 1 \bmod p$,
- by *tame relations* $r_v = (t_v)^{q_v}[t_v, s_v]$ with $q_v = |\mu(F_v)|$ $v \in S - S_0$ such that $v \equiv 1 \bmod p$,
- by *wild generators* $s_v, t_v, s_v^2, t_v^2, \ldots, s_v^{n_v}, t_v^{n_v}$ with $v \in S - S_0$ such that $v$ divides $p$,
- and by *wild relations* $r_v = (t_v)^{q_v}[t_v, s_v][t_v^2, s_v^2] \cdots [t_v^{n_v}, s_v^{n_v}]$ with $q_v = |\mu(F_v)| = p$, $v \in S - S_0$ such that $v$ divides $p$, $2n_v = [F_v : \mathbf{Q}_p] + 2$ (see [Se2] corollary 4.3).

The archimedean places do not appear because $p$ is odd.

**Remark 8.1.4** Wingberg's presentation yields a system of generators of $P_{G_S}$ but this system is not special. We have to work in order to obtain a 'good' system of generators, i.e for which we know enough information about the $A$-action.

8.2. *The* $\mathbf{Z}_p$*-cyclotomic extension of* $F$. Let us assume $\Gamma = \Gamma_{cyc} =< \hat{\gamma} >$. Recall that if $x \in P_{G_S}$ we denote $\tilde{x}$ its projection in $P_G$.

**Choice of the generators of** $P_G$. We have to choose a lift of $\hat{\gamma}$ in $P_G$. Recall that a tame generator $\tilde{s}_v$ comes from a generator of the local $\mathbf{Z}_p$-cyclotomic extension. We assume there exists $w \in S - S_0$ such that $N(w) \equiv 1 \bmod p$ and $N(w) \not\equiv 1 \bmod p^2$ hence $q_w = p$, so that we can assume that $\pi(\tilde{s}_w) = \hat{\gamma}$ generates the $\mathbf{Z}_p$-cyclotomic extension of $F$; if such a $w$ does not exists we can choose a wild generator $s_w$ such that $\pi(\tilde{s}_w) = \hat{\gamma}$. We denote $\gamma = s_w$ and $\tilde{\gamma} = \tilde{s}_w$.
We fix $\bar{\rho}(t_w^1) = \begin{pmatrix} 1 & 1 \\ 0 & 1 \end{pmatrix}$; this image depends upon the representation, but this choice does not change the method of computation.

In order to choose a 'good' minimal system of generators of $P_G$ let us express the image by the projection onto $\Gamma$ of the system of generators of $P_G$:
- The images by $\pi$ of $\tilde{t}_v$ or $\tilde{t}_v^i, \tilde{s}_v^i$, $2 \leq i \leq n_v$ are trivial by definition.
- To find the projection of the tame generator $\tilde{s}_v$ in $\Gamma$, it suffices to compute the number of $p^{th}$-roots of unity in the localization of $F$ in $v$; by Hensel's lemma it suffices to compute the number $q_v$ of $p^{th}$-roots of unity in the residue field. Hence the image of $\tilde{s}_v$ in $\Gamma$ is $\hat{\gamma}^{q'_v}$ and $q'_v$ has the same $p$-adic valuation as $q_v/p$.
- The projection of the wild generator $\tilde{s}_v$, in $\Gamma$ is $\hat{\gamma}^{q'_v}$.

Then $P_{G_S}$ admits the following 'good' system of generators (only the element $\gamma$ has a non

trivial image by $\pi$):

$$\{s_v,\ v \in S'\} \cup \{\gamma, t_w\} \cup$$

$$\{s'_v = \gamma^{-q'_v} s_v, t_v,\ v \in S - (S_0 \cup \{w\})\} \cup$$

$$\{s_v^i, t_v^i,\ v \in S - S_0,\ v|p,\ 2 \le i \le n_v\}$$

with the relations

$$r_w = (t_w)^p [t_w, \gamma]$$

$$r'_v = (t_v)^{q_v} [t_v, \gamma^{q'_v} s'_v],\ v \in S - S_0,\ v \not| p$$

$$r'_v = (t_v)^{q_v} [t_v, \gamma^{q'_v} s'_v][t_v^2, s_v^2] \cdots [t_v^{n_v}, s_v^{n_v}],\ v \in S - S_0, v|p$$

**Fox derivatives**. For $v \in S - \{w\}$ the Fox derivatives of a tame relation are

$$\frac{\partial r'_v}{\partial t_v} = \sum_{i=0}^{q_v} (t_v)^i - (t_v)^{q_v+1} \gamma^{q'_v} s'_v (t_v)^{-1}$$

$$\frac{\partial r'_v}{\partial \gamma} = \sum_{i=0}^{q'_v-1} (t_v)^{q_v+1} \gamma^i - \sum_{i=1}^{q'_v} (t_v)^{q_v+1} \gamma^{q'_v} s'_v (t_v)^{-1} (s'_v)^{-1} \gamma^{-i}$$

$$\frac{\partial r'_v}{\partial s'_v} = (t_v)^{q_v} t_v - (t_v)^{q_v} [t_v, s'_v \gamma^{q'_v}]$$

Since $\pi(\tilde{\gamma}) = \hat{\gamma}$ and all the other elements of the good system of generators of $P_G$ have a trivial image by $\pi$, the columns of the Fox matrix are of the form (for wild or tame relations):

$$(0, \ldots, 0, q_v - (\hat{\gamma}^{q'_v} - 1), 0, \ldots, 0),\ v \in S - S_0$$

**Action by** $A$. Let $l$ be a prime. To each $v|l$, $v \in S - S_0$, there corresponds a Demuškin group $\mathcal{P}_v$ one which we would like to determine the action of the decomposition group $A_l$. Let $q$ be the highest power of $p$ such that $F_v$ contains a primitive $q^{th}$-root of unity. Then $V_q = F_v^*/F_v^{*q}$ is a symplectic space relatively to the $q$-Hilbert symbol $< \cdot, \cdot >$, and $A_l$ acts on this space as a group of similarities, i.e

$$< a(x), a(y) >= \omega(a) < x, y >,\ \forall a \in A_l$$

where $\omega$ denotes the Teichmüller character.
Since the order of $A_l$ is prime to $p$, we have both tame ramification and semi-simplicity. Results of Borevič, Jakovlev and Koch ( [Ko2] proposition 6, [Ko2] proposition 9) on simplectic spaces with operators then allow us to decompose $V_q$ as a direct sum of 'hyperbolic planes' (in an obvious sense), $V_q = \oplus_\chi H_\chi$, where $H_\chi$ is generated by a hyperbolic pair $s_\chi, t_{\omega\chi^{-1}} \in V_{q,\chi}$.
Taking this into account and replacing the Demuškin refining process, is not enough to obtain a system of generators of $\mathcal{P}_v$ which verify the Demuškin relation (see [Bö2] proposition 2.5). We can also obtain such a nice relation in terms of generators with the desired action of $A_l$ via characters modulo $F_d^{(3,q)}$ where $F_d^{(i,q)}$ is defined by induction

$$F_d^{(1,q)} = F_d \text{ and } F_d^{(i+1,q)} = (F_d^{(i,q)})^q[F_d^{(i,q)}, F_d],$$

(see [Bö2] remark 2.4). For the wild generators, the action of $A$ is more ambiguous. We can easily be convinced of the difficulty of describing the action of $A$ upon recalling that in the Wingberg presentation, we forget a place $S_0$ over S.

**Image of the generators**. Let $\mu_p(F)$ and $\mu_p(F_v)$ be the sets of $p^{th}$-roots of unity in $F$ and $F_v$ respectively. Let $A'$ be the subgroup of $A$ of order two generated by $c$ and let $\tilde{\mathbf{F}}_p$ be the non trivial irreducible $\mathbf{F}_p[A']$-module.
Let $V^S = V^S_\emptyset$. Then we have $V^S \subset V^S_{S_0}$ and after [BoUl] proposition 3.2, if $V^S = 0$

$$\bar{P}_{G_S} \cong \mathrm{Ind}_{A'}^A \tilde{\mathbf{F}}_p \oplus \mathbf{F}_p \oplus \mathrm{Coker}(\mu_p(F) \to \oplus_{v\in S}\mu_p(F_v))$$

as an $\mathbf{F}_p[A]$-module. This description of $P_{G_S}$ gives us a special system of generators. In order to simplify our computation and to give a more explicit way to obtain the universal deformation ring we make the following assumption:
Let $S_*$ denotes the set of prime numbers such that $S$ is the set of places over $S_*$ in $F$. We assume that for all $l \in S_*$, $l$ doesn't split in $\mathbf{Q}(\mu_p)/\mathbf{Q}$. Thus we have

**Lemma 8.2.1** *Assume that* $\Gamma$ *contains the cyclotomic* $\mathbf{Z}_p$*-extension, that* $V^S_{S_0} = 0$ *and that*

*for all $l \in S_*$, $l$ doesn't split in $\mathbf{Q}(\mu_p)/\mathbf{Q}$, then*

$$X_{\infty,triv} = \{u \in X_\infty : \ a.u = u, \ \forall a \in A\} = 0$$

PROOF: Recall that a generator of the cyclotomic $\mathbf{Z}_p$-extension is invariant by complex conjugation.

Using [BoUl] proposition 3.2, in order to prove the lemma it suffices to show that the action of $A$ on $\mathrm{Ind}_{A'}^{A}\tilde{\mathbf{F}}_p$ and on $\mathrm{Coker}(\mu_p(F) \to \oplus_{v\in S}\mu_p(F_v))$ is not trivial. It is clear for the induced representation because the action of the complex conjugation $c$ is not trivial. It suffices to show that the following cokernel is zero:

$$\mathrm{Coker}(H^0(A, \mu_p(F)) \to H^0(A, \oplus_{v\in S}\mu_p(F_v)))$$

By definition $H^0(A, \mu_p(F)) = \mu_p(\mathbf{Q}) = 0$. After Shapiro's lemma, we have

$$H^0(A, \oplus_{v\in S}\mu_p(F_v)) \cong \oplus_{l\in S_*}H^0(A_l, \mu_p(\mathbf{Q}_l)) \cong \oplus_{l\in S_*}\mu_p(\mathbf{Q}_l)$$

Since $l \in S_*$ doesn't split in $\mathbf{Q}(\mu_p)/\mathbf{Q}$, $\mu_p(\mathbf{Q}_l) = 0$. ∎

Hence all the images by $\rho_G$ of the elements of $X_\infty$ are upper triangular. We denote:

$$\rho_G(\tilde{\gamma}) = \begin{pmatrix} 1+Y & 0 \\ 0 & 1+Y' \end{pmatrix}, \ \rho_G(\tilde{t}_w) = \begin{pmatrix} 1 & 1 \\ 0 & 1 \end{pmatrix}$$

$$\rho_G(\tilde{s}'_v) = \begin{pmatrix} 1 & Y'_v \\ 0 & 1 \end{pmatrix} \ v \in S \cup S' - \{w\}, \ \rho_G(\tilde{t}_v) = \begin{pmatrix} 1 & Y_v \\ 0 & 1 \end{pmatrix}, \ v \in S - \{w\}$$

and for wild generators with $v \in S - S_0$, $2 \le i \le n_v$

$$\rho_G(\tilde{s}^i_v) = \begin{pmatrix} 1 & {Y'}^i_v \\ 0 & 1 \end{pmatrix}, \ \rho_G(\tilde{t}^i_v) = \begin{pmatrix} 1 & Y^i_v \\ 0 & 1 \end{pmatrix}$$

possibly with some relations between the $Y_v, Y'_v, Y^i_v, {Y'}^i_v$ depending upon the action of $A$. To find the form of the ideal of relations $I$ we do not need to know precisely the action of $A$:

Using the images of the elements of $X_\infty$ we have linear relations between the $Y_v, Y'_v$ of the form

$$Y_v = f(\text{variables})$$

where $f$(variables) denotes a linear combinaison of the variables introduced by the images of the elements of $X_\infty$ (write $\rho_G(a \cdot x) = \rho_G(a)\rho_G(x)\rho_G(a)^{-1}$ $a \in A, x \in X_\infty$). Then the action of $A$ allows us to eliminate some variables. To summarize:

**Lemma 8.2.2** *In the Borel case, if* $\mathrm{Im}\bar{\rho}$ *is not diagonal, if* $V^S_{S_0} = 0$ *and if* $\Gamma = \Gamma_{cyc}$*, then a presentation by generators and relations of* $P_{G_S}$ *is known and the universal deformation ring of* $\bar{\rho}_G$ *reads*

$$R_G(\bar{\rho}_G) = \mathbf{Z}_p[[Y, Y', Y_v, Y'_{v'}, Y^i_{v''}, {Y'}^i_{v''}]]/I$$

*with* $v' \in S \cup S' - (S_0 \cup \{w\})$ *and* $v \in S - (S_0 \cup \{w\})$ *and* $v'' \in S - S_0$, $v''|p$, $2 \leq i \leq n_{v''}$*; where* $I$ *is the ideal of relations generated by*

$$\Big(q_v - \Big(\Big(\frac{1+Y}{1+Y'}\Big)^{q'_v} - 1\Big)\Big)Y_v,\ v \in S - (S_0 \cup \{w\})$$

$$\text{and } p - \Big(\frac{1+Y}{1+Y'} - 1\Big)$$

*and some linear combinations between the variables* $Y_v, Y_{v,}, Y^i_{v''}, {Y'}^i_{v''}$ *depending upon the action of* $A$.

Let replace $Y, Y'$ by $Y_1, Y'_1$ with $Y'_1 = \dfrac{1+Y'}{1+Y} - 1 = p$ which can be removed; then we have

$$R_G(\bar{\rho}_G) \cong \mathbf{Z}_p[[Y_1, \ldots, Y_{d''}]]/I'$$

where $d'' = \dim_{\mathbf{F}_p} H^1(P_G, Ad\bar{\rho}_G)$ and $I'$ is an ideal generated by relations of the form $pY_i$ with $r \geq 1$.

**Remark 8.2.3** Let $\bar{R}_G(\bar{\rho}_G) = R_G(\bar{\rho}_G) \otimes \mathbf{F}_p$. Then $\bar{R}_G(\bar{\rho}_G)$ is free. It shows that some deformations of the residual representation $\bar{\rho}$ cannot lift to characteristic zero rings, while others can (see [Ti] subsection 5.2).

8.3. *The maximal free pro-p-extension of $F$.* If $P_{G_S}$ admits a Wingberg presentation, is not needed $w_\Gamma = v_\Gamma$ because then the relations are simple, though the action by $A$ on the generators introduced by Wingberg's presentation is not easy to write.
Let $\Gamma = \Gamma_{max}$.
We use the presentation of $P_{G_S}$ and the notations of paragraph 8.1. Then we can choose $\Gamma_{max}$ such that

$$(\pi(\tilde{s}_v)\ v \in S' \cup S - S_0,\ \pi(\tilde{s}^i_{v'})\ v' \in S - S_0\ v'|p\ 2 \leq i \leq n_v)$$

is a minimal system of generators of $\Gamma$ (see lemma 4.7 [Ya1]).

**Fox derivatives**. The computation of the Fox derivatives as above gives the columns of the Fox matrix for a wild relation:

$$(0, \ldots, 0, q_v - (\pi(\tilde{s}_v) - 1), -(\pi(\tilde{s}^2_v) - 1), \ldots, -(\pi(\tilde{s}^{n_v}_v) - 1), 0, \ldots, 0), v \in S - S_0.$$

**Action by** $A$. Using the previous discussion on the action by $A$ we have a good description modulo $F_d^{(p,3)}$ (with the notations of 8.1.) of $P_{G_S}$ by generators and relations with action by $A$:

• Tame generators. Let $l$ be a prime $l \neq p$. We assume that $S$ contains all places over $l$ in $F$ and that $A/A_l \cong \mathrm{Gal}(\mathbf{Q}(\zeta_p)/\mathbf{Q})$. Then $A \cong A_l \times \mathbf{F}_p^*$.
Let $v$ be a place over $l$ in $F$. The local action by $A_l$ on $\mathcal{P}_v \cong < \hat{t}_v, \hat{s}_v | \hat{r}_v >$ is known and $A$ acts on $\oplus_{v'|l}\mathcal{P}_{v'}$ by the induced representation $\mathrm{Ind}_A^{A_l}\mathcal{P}_v^{ab}$. We can lift the local generators to $s_v, t_v$ keeping the action by $A_l$ modulo $F_d^{(2,p)}$; hence modulo $F_d^{(2,p)}$, $A$ permutes the $t_v$, $v|l$ (resp. $s_v$).

• Wild generators. We assume there are at most two places over $p$ in $S_F$, i.e we have at most one wild relation. In the same way as for the tame generators, we can lift modulo $F_d^{(2,p)}$ the wild generators with the action by $A$ described in the previous subsection: wild generators $t_v, s_v, t_v^2, s_v^2, \ldots, t^{n_v}, s_v^{n_v}$ where $t_v \in P_{G,\omega}$, $s_v \in P_{G,triv}$, $t_v^i \in P_{G,\omega\chi_i^{-1}}$, $s_v^i \in P_{G,\chi_i}$ modulo $F^{(2,p)}$ (you get image by $\rho_G$ modulo $\mathrm{Gl}_2(R)^{(2,p)}$) with a relation $t_v^p[t_v, s_v][t_v^2, s_v^2]\cdots[t_v^{n_v}, s_v^{n_v}]$ modulo $F_d^{(3,p)}$.

Knowing the Fox derivatives and the action of $A$ we can now compute an approximation of the universal deformation ring $R_G(\bar{\rho})$ (for a discussion of this approximation see [Bö2]).

*Acknowledgement*: The author wishes to thank Professor Gillard for many valuable discussions and Professor Nguyen Quang Do for his constant help and enjoyable advice. The author also acknowledges G. Böckle's critical comments on the manuscript.

*References*

[Bö1] G. Böckle. Explicit Universal Deformations of Even Galois Representations. Preprint (1996).

[Bö2] G. Böckle. Demušhkin groups with group actions and applications to deformations of local Galois representations. Preprint (1997).

[Bo1] N. Boston. Explicit deformation of Galois representations. *Invent. math.* **103** (1991), 181-196.

[Bo2] N. Boston. Deformations of Galois representation associated to the cusp form Δ. *Séminaire de Théorie des Nombres de Paris* (1987), 51-62.

[Bo3] N. Boston. Deformation Theory of Galois representations. Ph.D. thesis. Harvard University Cambridge, Massachusetts (1987).

[Bo4] N. Boston. Families of Galois representations. Increasing the ramification. *Duke Mathematical Journal* **66** (1992), 357-367.

[BoUl] N.Boston and S.V Ullom. Representations related to CM elliptic curves. *Math. Proc. Camb. Phil. Soc.* **113** (1993), 71-85.

[Br] A. Brumer. Pseudo-compact algebras, profinite groups and class-formations. *J. of Algebra* **4** (1966), 442-470.

[Ja] U. Jannsen. Iwasawa Modules up to isomorphism. *Adv. studies in Pure math.* **17** (1989), 171-207.

[Ko1] H. Koch. *Galoissche Theorie der p-Erweiterungen.* Springer-Verlag, New-York, Heidelberg, Berlin, 1970.

[Ko2] H. Koch. Über Darstellungsräume und die Struktur der multiplikativen Gruppe eines $p$-adischen Zahlkörpers. *Math. Nach.* **26** (1963), 67-100.

[Ma] B. Mazur. Deforming Galois representations. In *Galois groups over* **Q**. Y. Ihara, K. Ribet, J-P Serre eds., MSRI Publ. **16**, Springer-Verlag, New-York, Berlin, Heidelberg, 1987, 385-437.

[Ne] O. Neumann. On $p$-closed number fields and an analogue of Riemann's existence theorem. In *Algebraic Number Fiel.* Proc. Symp. London math. Soc., Univ. Durham 1975 (1977), 625-647.

[Ng1] T. Nguyen Quang Do. Formations de classes et modules d'Iwasawa. *Number Theory Noordwijkerhout 1983*, Lecture Notes in mathematics **1068** (1984), 167-185.

[Ng2] T. Nguyen Quang Do. Sur la structure galoisienne des corps locaux et la théorie d'Iwasawa. *Compositio Math.* **26** (1982), 85-119.

[Ng3] T. Nguyen Quang Do. Sur la $\mathbf{Z}_p$-torsion de certains modules galoisiens. *Ann. Inst. Fourier* **36** (1986), 27-46.

[Ng4] T. Nguyen Quang Do. Lois de réciprocité primitives. *Manuscripta Math.* **72** (1991), 307-324.

[Ra] R. Ramakrishna. On a variation of Mazur's deformation functor. *Compositio Math.* **87** (1993), 269-286.

[Sc] M. Schlessinger. Functors of Artin rings. *Trans. Am. Soc.* **130** (1968), 208-222.

[Se1] J.P Serre. Propriétés galoisiennes des points d'ordre fini des courbes elliptiques. *Inventiones math.* **15** (1972), 259-331.

[Se2] J.P Serre. Structure de certains pro-$p$-groupes. *Sém. Bourbaki* **252** (1963), 145-155.

[Se3] J.P Serre. Modular forms of weight one and Galois representations. In *Algebraic Number Fields.* Fröhlich eds., Acad. Press (1977), 193-268.

[Ti] J. Tilouine. *Deformations of Galois representations and Hecke algebras.* Metha Research Institute (1996).

[Wa] L.C Washington. *Introduction to cyclotomic Fields.* Springer Verlag **83**, New-York, Heidelberg, Berlin 1980.

[Wi] K. Wingberg. On Galois groups of $p$-closed algebraic number fields with restricted ramification II. *J. reine angew. Math.* **416** (1991), 187-194.

[Ya1] M. Yamagishi. A note on free pro-$p$-extensions of algebraic number fields. *Journal de Théorie des Nombres de Bordeaux* **5** (1993), 165-178.

[Ya2] M. Yamagishi. A note on free pro-$p$-extensions of algebraic number fields II. *Manuscripta Math.* **91** (1996), 231-233.